\documentclass[a4paper,10pt]{amsart}
\usepackage{amssymb}
\usepackage{amsthm}
\usepackage{zmacros}
\usepackage{url}
\usepackage[ps,dvips,2cell,matrix,arrow,line,tips]{xy} 
\usepackage{xspace}
\usepackage{srcltx}
\input{zmathmacros.cfg}
\input{abbreviations.cfg}
     
\author{Marco Zunino}
\title{Double Construction for Crossed Hopf Coalgebras}
\subjclass[2000]{81R50,16W30,57R56}
\date{Friday, 13 September 2002}
\address{Marco Zunino -- 
        Institut de Recherche Math\'ematiques Avanc\'ee (IRMA) -- 
        Universit\'e Louis Pasteur et CNRS -- 
        7,~rue Ren\'e Descartes, 67084 Strasbourg \textsc{cedex}, France. }
\email{zunino@math.u-strasbg.fr}

\begin{document}

\bibliographystyle{amsplain}

\begin{abstract}
 We provide an analog of the Drinfeld quantum double construction
 in the context of crossed Hopf group coalgebras introduced by Turaev.
 We prove that, provided the base group is finite, the double of a 
 semisimple crossed Hopf group coalgebra is both modular and 
 unimodular.
\end{abstract}

\maketitle

\tableofcontents

\section*{Introduction}
Recently, Turaev~\cite{Tur-pi,Tur-CPC} (see also Le and 
Turaev~\cite{LeTur} and Virelizier~\cite{Virelizi2}) generalized 
the notion of a TQFT and Reshetikhin-Turaev invariants to the case
of $3$\nobreakdash-\hspace{0pt}manifolds endowed with 
a homotopy classes of maps to $K(\pi,1)$, where $\pi$ is a group.
One of the key points in~\cite{Tur-pi} is the notion of a 
\textit{crossed Hopf \picoalg,} here called a 
\textit{Turaev coalgebra} or, briefly, a \textit{\Tcoalg\/} 
(see Section~\ref{sec:CGC}). 
As one can use categories of representations of modular Hopf algebras 
to construct Reshetikhin-Turaev 
invariants of $3$\nobreakdash-\hspace{0pt}manifolds, 
one can use categories of representations of modular \Tcoalgs to construct 
homotopy invariants of maps from $3$\nobreakdash-\hspace{0pt}manifolds to
$K(\pi,1)$. 
Similarly, to construct Virelizier Hennings-like homotopy 
invariants~\cite{Virelizi2}, we need a ribbon \Tcoalg $H$ such that 
the neutral component $H_{1}$ is  unimodular~\cite{Schneider}.

Roughly speaking, a \Tcoalg $H$ is a 
family $\{H_{\alpha}\}_{\alpha\in\pi}$ of algebras endowed with a \textit{comultiplication}
$\Delta_{\alpha,\beta}\colon H_{\alpha\beta}\to H_{\alpha}\otimes H_{\beta}$, a \textit{counit} $\varepsilon\colon\Bbbk\to H_{1}$ 
(where $1$ is the neutral element of $\pi\/$), and an \textit{antipode}
$s_{\alpha}\colon H_{\alpha}\to H_{\alpha^{-1}}$. 
It is also required that $H$ is endowed with a family 
of algebra isomorphisms $\varphi^{\alpha}_{\beta}=\varphi_{\beta}\colon H_{\alpha}\to H_{\beta\alpha\beta^{-1}}$, 
the \textit{conjugation,} compatible with the above structures 
and such that $\varphi_{\beta\gamma}=\varphi_{\beta}\circ\varphi_{\gamma}$.
When $\pi=\{1\}$, we recover the definition of a Hopf algebra. 
A \Tcoalg $H$ is \textit{of finite type} when every $H_{\alpha}$ 
is finite-dimensional. $H$ is \textit{totally-finite} 
when the direct sum $\bigoplus_{\alpha\in\pi}H_{\alpha}$ is finite-dimensional.
A \textit{universal \Rmatrix}
and a \textit{twist} for a \Tcoalg $H$ are, respectively, families
$R=\bigl\{\xi_{(\alpha)}\otimes\zeta_{(\beta)}=R_{\alpha,\beta}\in H_{\alpha}\otimes H_{\beta}\bigr\}_{\alpha,\beta\in\pi}$ and
$\theta=\{\theta_{\alpha}\in H_{\alpha}\}_{\alpha\in\pi}$ satisfying axioms that explicitly involve 
the conjugation (see, respectively, Sections~\ref{fur-Alina} 
and~\ref{fur-Alina-2}). $H$ is \textit{modular} when it is ribbon and
its component $H_1$ is modular (see Section~\ref{fur-Alina-3}).
Properties of \Tcoalgs are studied in~\cite{Virelizi}.

 Starting from a finite-dimensional Hopf algebra $H$, Drinfeld~\cite{Drn}
 showed how to obtain a quasitriangular Hopf algebra $D(H)$, the 
 \textit{quantum double of $H$,} such that the following conditions 
 are satisfied (for details, see, e.g.,~\cite{Kas}).
 \begin{itemize}
  \item There are embeddings of Hopf algebras $i\colon H\hookrightarrow D(H)$
        and $j\colon H^{\ast\cop}\hookrightarrow D(H)$.
  \item The linear map
        $H\otimes H^{\ast\cop}\xrightarrow{i \otimes j} D(H)\otimes D(H)\xrightarrow{\mu} D(H)$
        is bijective, where $\mu$ is the multiplication in $D(H)$.
  \item The universal \Rmatrix of $D(H)$ is the image of the
        canonical element of $H\otimes H^{\ast\cop}$ under the
        embedding $i\otimes j\colon H\otimes H^{\ast\cop}\hookrightarrow D(H)\otimes D(H)$.
 \end{itemize}

In this paper we provide an analog of the quantum double for \Tcoalg $H$. 
We start (Section~\ref{sec:QD}, Theorem~\ref{thm:Asterix}) with an abstract
definition of $D(H)$, as a solution of a  universal problem,
analogous to the above definition of the standard quantum double.
The notion of a \Tcoalg in not 
self-dual, i.e., given a \Tcoalg $H=\{H_{\alpha}\}_{\alpha\in\pi}$, the family 
$\Houter_{\alpha}=\{H^{\ast}_{\alpha}\}_{\alpha\in\pi}$ does not have a natural structure of a \Tcoalg. 
However, when $H$ is of finite type, we introduce a \Tcoalg $\Hinner$, 
the \textit{inner dual of $H$,} 
that, in many aspects, in particular in the construction of $D(H)$,
plays the role of a dual for $H$
(see Section~\ref{sec:duals}). The components of $\Hinner$ are all
isomorphic as algebras. In particular, as a vector space,
$\Hinner_{\alpha}=\bigoplus_{\beta\in\pi}H^{\ast}_{\beta}$ for any $\alpha\in\pi$.
 
An explicit description of $D(H)$
is given in Theorems~\ref{thm:dual} and~\ref{thm:DH}.
In particular, every component $D_{\alpha}(H)$ of $D(H)$ is,
as a vector space, $H_{\alpha^{-1}}\otimes\bigoplus_{\beta\in\pi}H^{\ast}_{\beta}$. So,
in general, $D(H)$ is of finite type if and only if 
$H$ is totally-finite and, in that case, $D(H)$ is also totally-finite.
The product in $D_{\alpha}(H)$ is obtained by setting
         \begin{gather*}
           (1_{\alpha^{-1}}\circledast f)\,(h\circledast\varepsilon)  = h\circledast f\qquad\text{and}\\
           (h\circledast\varepsilon)\, (1_{\alpha^{-1}}\circledast f) =
           h''_{(\alpha^{-1})}\circledast \bigl\langle f, 
           s_{\gamma}^{-1}(h'''_{(\gamma^{-1})})\_\varphi_{\alpha}(h'_{(\alpha^{-1}\gamma\alpha)})\bigr\rangle\text{.}
         \end{gather*}
for any $f\in\bigoplus_{\beta\in\pi}H^{\ast}_{\beta}$, and $h\in H_{\alpha^{-1}}$.
When $\pi=\{1\}$, we recover the standard definition of the quantum double 
of a Hopf algebra. We will see that both the product and the coproduct 
in $D(H)$ explicitly depend on the conjugation $\varphi$ of $H$.
 
 The quantum double $D(H)$ of a semisimple Hopf algebra $H$ over a field 
of characteristic $0$ is both semisimple~\cite{LR1,Radford} and 
modular~\cite{EG}. Here we prove that the double $D(H)$ of a semisimple 
\Tcoalg $H$ over a 
field of characteristic $0$ is semisimple if and only if $H$ is 
totally-finite, and, in that case, $D(H)$ is also modular 
(Theorem~\ref{thm:ssmod}). The key point in the proof 
is that, when $H$ is totally finite, it gives rise to a graded Hopf
algebra $H_{\pk}=\bigoplus_{\alpha\in\pi}H_{\alpha}$, the \textit{packed form of $H$.}
The Hopf algebras $D(H_{\pk})$ and $\bigl(D(H)\bigr)_{\pk}$ are different,
but it is always possible to embed $D_{1}(H)$ in $D(H_{\pk})$ as an algebra. 
 
 Last section is devoted to the generalization of the Reshetikhin-Turaev 
ribbon construction. Given a quasitriangular Hopf algebra $H$, Reshetikhin and Turaev~\cite{RT} embedded it into a ribbon Hopf algebra $\RT(H)$ 
that is a quotient of the polynomial algebra $H[\theta]$.
Starting from a quasitriangular \Tcoalg $H$ 
(not necessarily of finite type), we construct a ribbon \Tcoalg $\RT(H)$ 
such that, when $\pi=\{1\}$, we recover the Reshetikhin-Turaev construction
for Hopf algebras. As a corollary, for any \Tcoalg $H$ of finite type
(but not necessarily totally-finite), 
we obtain a ribbon \Tcoalg $\RT\bigl(D(H)\bigr)$.

The quantum double and the ribbon extension have categorical 
counterparts~\cite{JS,KasTur,Street-double}. In~\cite{Zunino-center},
we provide a generalization of these constructions.
In particular, the generalization of the center construction
(corresponding to the quantum double in the algebraic case)
looks very natural and as the simplest way to adapt the standard
construction to the crossed case. Both the definition of $D(H)$ and
the definition of $RT(H)$ introduced here are coherent to the 
categorical constructions as we briefly explain.
In the standard case, starting from a tensor category
$\mathcal{C}$, we obtain a braided tensor category 
$\mathcal{Z}(\mathcal{C})$ via the center construction~\cite{JS}.
Now, take $\mathcal{C}=\Rep(H)$, the category of 
representations of a finite-dimensional Hopf algebra $H$. 
Both $\mathcal{Z}(\mathcal{C})$ and $\Rep\bigl(D(H)\bigr)$ 
are canonically isomorphic to the category of Yetter-Drinfeld modules
over $H$ (see~\cite{Majid-YD,Yetter}). The corresponding result
is still true when we consider a \Tcoalg $H$ instead of an Hopf
algebra, that is, also in the crossed case, we have
a commutative diagram 
 \begin{equation*}
 \vcenter{\xymatrix{\Rep(H)\ar[r]^{\mathcal{Z}\ \ \ \ } 
           & \mathcal{Z}\bigl(\Rep(H)\bigr)\ar@{=}[r]
           & \YD\,(H)\ar@{=}[r] 
           & \Rep\bigl(D(H)\bigr)\\ 
          H\ar[u]^{\Rep}\ar[rrr]_{D} & \ & \ & D(H)\ar[u]_{\Rep}}}
\end{equation*}
(where, of course, we also need to generalize the notions of
a tensor category and of a Yetter-Drinfeld module).
Similar results are obtained for the ribbon construction.\medskip

\paragraph{\scshape Acknowledgments} 
The author wants to thank his adviser V.\@~Turaev for his
stimulating and constructive direction in the research. 
The author also wants to thank A.\@~Brugui\`eres, 
B.\@~Enriquez, Ch.\@~Kassel, H.J.\@~Schneider, L.\@~Vainerman,
and A.\@~Virelizier for their useful remarks.
  
The author was partially supported by the 
\textsc{INdAM, Istituto Nazionale di Alta Matematica,} 
Rome.

\section{\protect\Tcoalgs}\label{sec:CGC}
 We recall the notion of a 
\Tcoalg as in~\cite{Tur-CPC}. \Tcoalgs are a generalization 
of standard notion Hopf algebras (for a general reference to
Hopf algebras, see~\cite{Sweedler,Abe} or, for a 
modern introduction,~\cite{Schneider}).
A generalization of the Heynemann\hspace{0pt}-\hspace{0pt}Sweedler 
notation~\cite{Sweedler} is also provided.
At the end of this section, we study some 
properties of the antipode of a \Tcoalg.
  
\paragraph{\scshape Basic definitions}
 Let $\Bbbk$ be a commutative field and let $\pi$ be a group. 
 A \textit{\Tcoalg $H$ \textup{(}over $\pi$ and $\Bbbk$\/\textup{)}} 
 is given by the following data.
 \begin{itemize}
  \item For any $\alpha\in\pi$, an associative \kalg $H_{\alpha}$, called the
        \textit{\alphath component of $H$.} The multiplication 
        is denoted $\mu_{\alpha}\colon H_{\alpha}\otimes H_{\alpha}\to H_{\alpha}$
        and the unit is denoted 
        $\eta_{\alpha}\colon \Bbbk\to H_{\alpha}$, with $1_{\alpha}\eqdef\eta_{\alpha}(1)$. 
  \item A family of algebra morphisms 
        $\Delta=\bigl\{\Delta_{\alpha,\beta}\colon H_{\alpha\beta}\to H_{\alpha}\otimes H_{\beta}\bigr\}_{\alpha,\beta\in\pi}$,
        called \textit{comultiplication,} that is \textit{coassociative} 
        in the sense that, for any $\alpha,\beta,\gamma\in\pi$, we have
        $(H_{\alpha} \otimes \Delta_{\beta, \gamma})\circ\Delta_{\alpha,\beta\gamma}=(\Delta_{\alpha,\beta}\otimes H_{\gamma})\circ\Delta_{\alpha\beta,\gamma}$.
  \item An algebra morphism $\varepsilon\colon H_{1}\to  \Bbbk$,
        called \textit{counit,} such that, for any $\alpha\in\pi$, we have
        $(\varepsilon \otimes H_{\alpha})\circ\Delta_{1,\alpha}=\Id$ and $(H_{\alpha} \otimes \varepsilon)\circ\Delta_{\alpha,1}=\Id$.
 \item A set of algebra isomorphisms 
       $\varphi= \bigl\{\varphi_{\beta}^{\alpha}\colon H_{\alpha}\to H_{\beta\alpha\beta^ {-1}} \bigr\}_{\alpha,\beta\in\pi}$
       called \textit{conjugation.}
       When not strictly necessary, the upper index will be omitted. 
       We require that $\varphi$ satisfies the following conditions. 
       \begin{itemize}
         \item $\varphi$ is \textit{multiplicative,} i.e., 
               $\varphi_{\beta}\circ\varphi_{\gamma}=\varphi_{\beta\gamma}$, for any $\beta,\gamma\in\pi$.
               It follows that  $\varphi_{1}^{\alpha}=\Id$, for any $\alpha\in\pi$.
         \item $\varphi$ is \textit{compatible with $\Delta$,} 
               i.e, $\Delta_{\gamma\alpha\gamma^{-1},\gamma\beta\gamma^{-1}}\circ\varphi_{\gamma}=(\varphi_{\gamma}\otimes\varphi_{\gamma})\circ \Delta_{\alpha,\beta}$,
               for any $\alpha,\beta,\gamma\in\pi$.
         \item $\varphi$ is \textit{compatible with $\varepsilon$,} i.e., 
               $\varepsilon\circ\varphi_{\gamma}=\varepsilon$ for any $\gamma\in\pi$.
       \end{itemize}
  \item Finally, a set of \klin morphisms  
        $s=\{s_{\alpha} \colon H_{\alpha}\to H_{\alpha^{-1}}\}_{\alpha\in\pi}$,
        the \textit{antipode,} such that, for any $\alpha\in\pi$,
        \begin{equation}\label{e:antipode}
         \mu_{\alpha}\circ(s_{\alpha^{-1}}\otimes H_{\alpha})\circ \Delta_{\alpha^{-1},\alpha}= \eta_{\alpha}\circ\varepsilon =
         \mu_{\alpha}\circ(H_{\alpha}\otimes s_{\alpha^{-1}})\circ\Delta_{\alpha,\alpha^{-1}}\text{.}
        \end{equation}
 \end{itemize}
 
 \noindent (The compatibility of the antipode with the conjugation 
 isomorphisms is a consequence of the axioms.)
 In~\cite{Tur-CPC}, a \Tcoalg is called 
 a \textit{crossed group Hopf coalgebra.}

 We{\label{pag:tf}} say that $H$ is \textit{of finite-type} when any 
 component $H_{\alpha}$ (with $\alpha\in\pi\/$) is a finite-dimensional \kvector space. 
 We say that $H$ is \textit{totally-finite} when $\dim_{\Bbbk}\bigoplus_{\alpha\in\pi}H_{\alpha}<\infty$, 
 i.e., when $H$ is of finite-type and almost all the $H_{\alpha}$ are zero. 
 It was proved in~\cite{Virelizi} that the antipode of a finite-type 
 \Tcoalg is always bijective.

 We observe that the component $H_{1}$ of a \Tcoalg $H$ is a Hopf algebra
 in the usual sense. We also observe that, for $\pi=\{1\}$, we recover the 
 usual notion of a Hopf algebra.

 \begin{ex}[\THcoalgs]\label{ex:quasiconstant}
  Let $H_{1}$ be a Hopf algebra, with comultiplication $\Delta_{1}$, 
  counit $\varepsilon_{1}$, and antipode $s_{1}$, endowed with a group morphism
   \begin{equation*}
       \map{\varphi^{1}}{\pi}{\Aut(H_{1})}{\alpha}{\varphi^{1}_{\alpha}}\text{,}
   \end{equation*}
 where $\Aut(H_{1})$ is the group of Hopf algebra
 automorphisms of $H_{1}$. We obtain a \Tcoalg $H$ by setting
 (for any $\alpha,\beta\in\pi\/$) $H_{\alpha}\eqdef H_{1}$ as algebra, $\Delta_{\alpha,\beta}\eqdef\Delta_{1}$,  
 $\varepsilon\eqdef\varepsilon_{1}$, $s_{\alpha}=s_{1}$, and $\varphi^{\alpha}_{\beta}=\varphi^{1}_{\beta}\colon H_{\alpha}\to H_{\beta\alpha\beta^{-1}}$.
 The \Tcoalg $H$ is called the \textit{\THcoalg based on $H_{1}$}
 (the H stand for Hopf, since $H$ is also an Hopf algebra).
\end{ex}

\paragraph{\scshape Coopposite \Tcoalg}
 Let $H$ be a \Tcoalg with invertible antipode. 
 The \textit{coopposite \Tcoalg}$H^{\cop}$ 
 is the \Tcoalg defined as follows.
      \begin{itemize}
       \item For any $\alpha\in\pi$, we set $H^{\cop}_{\alpha}\eqdef H_{\alpha^{-1}}$ 
             as an algebra.
       \item The comultiplication $\Delta^{\cop}$ is obtained by setting, 
             for any $\alpha,\beta\in\pi$,
             \begin{equation*}
              \Delta^{\cop}_{\alpha,\beta}\eqdef\Bigl(H^{\cop}_{\alpha\beta}=H_{\beta^{-1}\alpha^{-1}}
              \xrightarrow{\Delta_{\beta^{-1},\alpha^{-1}}}H_{\beta^{-1}}\otimes H_{\alpha^{-1}} 
              \xrightarrow{\sigma}H_{\alpha^{-1}}\otimes H_{\beta^{-1}}=
              H^{\cop}_{\alpha}\otimes H^{\cop}_{\beta}\Bigr)
              \text{.}
             \end{equation*}
             $\Delta^{\cop}$ has a counit $\varepsilon^{\cop}=\varepsilon$.
       \item The antipode $s^{\cop}$ is obtained by setting
             $s^{\cop}_{\alpha}=s^{-1}_{\alpha}\colon$, for any $\alpha\in\pi$.
       \item The conjugation  $\varphi^{\cop}$ is obtained by setting
             $\varphi^{\cop}_{\beta}\eqdef\varphi_{\beta}$, for any $\beta\in\pi$.
      \end{itemize}
      In particular, when $H$ is a \THcoalg, $H^{\cop}$ is
      the \THcoalg based on $H_1^{\cop}$.

\paragraph{\scshape Heynemann-Sweedler notation}
 The coassociativity of $H$ allows us to introduce 
 an analog of the Heynemann\hspace{0pt}-\hspace{0pt}Sweedler 
 notation, see~\cite{Sweedler}. Given $\alpha_{1},\dots, \alpha_{n}\in\pi$ and defining
 \begin{equation*}\begin{split}
  \Delta_{\alpha_{1},\alpha_{2},\dots,\alpha_{n}}\eqdef \Bigl(& H_{\alpha_{1}\alpha_{2}\cdots\alpha_{n}}
  \xrightarrow{\Delta_{\alpha_{1},\alpha_{2}\cdots\alpha_{n}}} H_{\alpha_{1}}\otimes H_{\alpha_{2}\cdots\alpha_{n}}
    \xrightarrow{H_{1}\otimes\Delta_{\alpha_{2},\alpha_{3}\cdots\alpha_{n}}}\\
    & H_{1}\otimes H_{2}\otimes H_{\alpha_{3}\cdots\alpha_{n}}
  \to\ \cdots\ \to H_{\alpha_{1}}\otimes H_{\alpha_{2}}\otimes \cdots\otimes H_{\alpha_{n}}\Bigr)\text{,}
 \end{split}\end{equation*}
 for any $h\in H_{\alpha_{1}\alpha_{2}\cdots\alpha_{n}}$ we set
 \begin{equation*}
  h'_{(\alpha_{1})}\otimes h''_{(\alpha_{2})}\otimes\cdots\otimes {h'}{}^{n}_{(\alpha_{n})}\eqdef
  \Delta_{\alpha_{1},\alpha_{2},\dots,\alpha_{n}}(h)\text{.}
 \end{equation*}
 Let $M$ be a vector space over $\Bbbk$ and suppose that 
 $f\colon H_{\alpha_{1}}\times H_{\alpha_{2}}\times\cdots\times H_{\alpha_{n}}\to M$
 is a \kmlin map. Let $\hat{f}$ denote the linear map
 $H_{\alpha_1}\otimes H_{\alpha_2}\otimes\cdots\otimes H_{\alpha_n}\to M$ induced by $f$. 
 We introduce the notation
 \begin{equation*}
  f\Bigl(h'_{(\alpha_{1})}, h''_{(\alpha_{1})},\cdots, {h'}{}^{n}_{(\alpha_{n})}\Bigr)\eqdef
  \hat{f}\bigl(\Delta_{\alpha_{1},\alpha_{2},\dots,\alpha_{n}}(h)\bigr)\text{.}
 \end{equation*}
 For simplicity, we also suppress the subscript $(\alpha_{i})$ when $\alpha_{i}=1$.
   
 \paragraph{\scshape Properties of the antipode}\label{convolution}
 Let $H$ be a \Tcoalg and 
 let $A$ be an algebra with multiplication $\mu_{A}$ and 
 unit $\eta_{A}$. We define a \textit{convolution algebra $\Conv(H,A)$} 
 (see~\cite{Virelizi}) in the following way. 
 As a vector space,
 $\Conv(H,A)\eqdef\bigoplus_{\beta\in\pi}\Hom_{\Bbbk}(H_{\beta},A)$.
 The multiplication in $\Conv(H,A)$ is obtained by setting, 
 for any $\beta_{1},\beta_{2}\in\pi$, $f_{1}\in\Hom_{\Bbbk}(H_{\beta_{1}},A)$, 
 $f_{2}\in\Hom_{\Bbbk}(H_{\beta_{2}}, A)$, and $h\in H_{\beta_{1}\beta_{2}}$, 
 \begin{equation*}
  (f_{1}\ast f_{2})(h)=
  f_{1}\bigl(h'_{(\beta_1)}\bigr)f_{2}\bigl(h''_{(\beta_2)}\bigr)\text{.}
 \end{equation*}
 With this multiplication, $\Conv(H,A)$ becomes 
 an associative algebra.
 
 For any $\alpha\in\pi$, we introduce the notation 
 $\Conv_{\alpha}(H)\eqdef \Conv(H,H_{\alpha})$.
 It is clear that~\eqref{e:antipode}
 is equivalent to say that $s_{\alpha^{-1}}$ is a two-sided inverse 
 of the identity morphism of $H_{\alpha}$ in the convolution algebra
 $\Conv_{\alpha}(H)$. Thus, we can reformulate the axiom for the antipode
 of $H$ by requiring that, for any $\alpha\in\pi$, the identity morphism of
 $H_{\alpha}$ is invertible in $\Conv_{\alpha}(H)$. This also proves that the
 antipode of a \Tcoalg is unique.
 
 \begin{lemma}
  Let $H$ be a \Tcoalg. The antipode $s$ of $H$ is both 
  antimultiplicative and anticomultiplicative, i.e.,
  \begin{equation}\label{e:s-antim}
   s_{\alpha}(hk) = s_{\alpha}(k)s_{\alpha}(h)
  \end{equation}
  for any $\alpha\in\pi$, $h,k\in H_{\alpha}$, and
  \begin{equation}\label{e:s-anticm}
   s_{\alpha}\bigl(h'_{(\alpha)}\bigr)\otimes s_{\alpha}\bigl(h''_{(\beta)}\bigr)=
   \bigl(s_{\alpha\beta}(h)\bigr)''_{(\alpha^{-1})}\otimes\bigl(s_{\alpha\beta}(h)\bigr)'_{(\beta^{-1})}
  \end{equation}
  for any $h\in H_{\alpha\beta}$.
 \end{lemma}
 
The proof can be obtained as in the standard case (i.e., for $\pi=\{1\}\/$).
 
\section{\protect\Talgs}
 When we dualize the axioms of a \Tcoalg, we obtain 
 the notion of a \Talg.
 A \Talg $H$ can be equivalently described as a Hopf algebra 
 endowed with a family of automorphisms, the 
 \textit{packed form of $H$.} In this section, we present both 
 approaches.
 
\paragraph{\scshape Basic definitions}
 A \Talg $H$ is a family 
 $\bigl\{(H_{\alpha},\Delta_{\alpha},\eta_{\alpha})\bigr\}_{\alpha\in\pi}$ of \kcoalgs, 
 endowed with the following data.
 \begin{itemize}
  \item A family of coalgebra morphisms
        $\mu=\{ \mu_{\alpha,\beta}\colon H_{\alpha}\otimes H_{\beta}\to  H_{\alpha\beta}\}_{\alpha,\beta\in\pi}$,
        called \textit{multiplication,} that is 
        \textit{associative,} in the sense that, for any $\alpha,\beta,\gamma\in\pi$,
        \begin{equation}\label{e:associativity}
         \mu_{\alpha\beta,\gamma}\circ (\mu_{\alpha,\beta}\otimes H_{\gamma})=
         \mu_{\alpha,\beta\gamma}\circ (H_{\alpha}\otimes\mu_{\beta,\gamma})\text{.}
        \end{equation}
        Given $h\in H_{\alpha}$ and $k\in H_{\beta}$, with $\alpha,\beta\in\pi$, we set
        $hk\eqdef\mu_{\alpha,\beta}(h,k)$.
        With this notation,~\eqref{e:associativity} 
        can be simply rewritten as
        $(hk)l = h(kl)$
        for any $h\in H_{\alpha}$, $k\in H_{\beta}$, $l\in H_{\gamma}$ and $\alpha,\beta,\gamma\in\pi$. 
  \item An algebra morphism
        $\eta\colon \Bbbk\to  H_{1}$,
        called \textit{unit,} such that, if we set 
        $1\eqdef\eta(1_{\Bbbk})$, then, for any $h\in H_{\alpha}$ 
        (with $\alpha\in\pi\/$), we have $1h = h = h1$.
 \item A set of coalgebra isomorphism 
       $\psi= \bigl\{ \psi_{\beta}^{\alpha}\colon H_{\alpha}\to H_{\beta\alpha\beta^ {-1}} \bigr\}_{\alpha,\beta\in\pi}$,
       called \textit{conjugation.} Also in this case, 
       when not strictly necessary, the upper index will be omitted. 
       We require that $\psi$ satisfies the following conditions.
       \begin{itemize}
        \item  $\psi$ is \textit{multiplicative,} i.e., 
               for any $\alpha$, $\beta$, and $\gamma\in\pi$, we have
               $\psi_{\beta}\circ\psi_{\gamma}=\psi_{\beta\gamma}$.
              It follows that, for any $\alpha\in\pi$, we have $\psi_{1}^{\alpha}=\Id$.
       \item  $\psi$ is \textit{compatible with $\mu$,} i.e,  for any $\beta\in\pi$,
               we have $\psi_{\beta}(hk) = \psi_{\beta}(h)\psi_{\beta}(k)$.
        \item $\psi$ is \textit{compatible with $\eta$,} i.e., for any $\beta\in\pi$,
              we have $\psi_{\beta}(1) = 1$.
       \end{itemize}
 \item Finally, a set of linear isomorphisms  
       $S=\{S_{\alpha} \colon H_{\alpha}\to H_{\alpha^{-1}}\}_{\alpha\in\pi}$,
       the \textit{antipode,} such that, for any $\alpha\in\pi$, we have
       $\mu_{\alpha^{-1},\alpha}\circ (S_{\alpha}\otimes H_{\alpha})\circ \Delta_{\alpha}=
        \eta\circ\varepsilon_{\alpha}=\mu_{\alpha,\alpha^{-1}}\circ (H_{\alpha}\otimes S_{\alpha})\circ \Delta_{\alpha}$.
 \end{itemize}
 
\noindent For any $\alpha\in\pi$, the coalgebra $H_{\alpha}$ is called 
the \textit{\alphath component of $H$.}

\paragraph{\scshape Packed form of a \protect\Talg}\label{par:packed-alg}
 Let $H$ be a \Talg. We define a Hopf algebra 
 $H_{\pk}$, that we call the \textit{packed form of $H$,} as follows.
 As a coalgebra, $H_{\pk}$ is the direct sum of the components of $H$. 
 For any $\alpha\in\pi$, we denote by $i_{\alpha}$ the inclusion of 
 $H_{\alpha}$ in $H_{\pk}$.
 The multiplication $\mu_{\pk}$ of $H_{\pk}$ is the colimit, in the 
 category of vector spaces, 
 $\mbox{$\varinjlim$}_{\alpha,\beta\in\pi}(i_{\alpha\beta}\circ\mu_{\alpha,\beta})$, 
 i.e., the only \klin map from $H_{\pk}\otimes H_{\pk}$ to $H_{\pk}$ 
 such that the restriction on $H_{\alpha}\otimes H_{\beta}\subset H_{\pk}\otimes H_{\pk}$ 
 coincides with $\mu_{\alpha,\beta}$,
 \begin{equation*}
  \vcenter{\xymatrix@C=1.8pc{ H_{\alpha}\otimes H_{\beta}\ar@{^{(}->}[r]\ar[d]_(.55){\mu_{\alpha,\beta}} &
  {\displaystyle\bigoplus_{\gamma,\delta\in\pi}(H_{\gamma}\otimes H_{\delta})=\Biggl(\bigoplus_{\gamma\in\pi}H_{\gamma}\Biggr)\otimes 
  \Biggl(\bigoplus_{\delta\in\pi}H_{\delta}\Biggr)= H_{\pk}\otimes H_{\pk}}\ar[d]^(.55){\mu_{\pk}}\ \ \\
  H_{\alpha\beta}\ar@{^{(}->}[r] & H_{\pk}}}\!\!\!\text{.}%))
 \end{equation*}
 The unit $\eta_{\pk}$ of $H_{\pk}$ is given by $i_{1}\circ\eta_{H}$, i.e., 
 $1_{\pk}=1\in H_{1}\subset H_{\pk}$. Finally, the antipode $S_{\pk}$ of $H_{\pk}$ 
 is given by the sum $\sum_{\alpha\in\pi}S_{\alpha}$. The Hopf algebra $H_{\pk}$ is 
 endowed with a group morphism $\map{\psi}{\pi}{\Aut(H_{\pk})}{\alpha}{\psi_{\pk,\alpha}}$, 
 where $\psi_{\pk,\alpha}\eqdef\sum_{\beta\in\pi}\psi_{\alpha}^{\beta}\colon\bigoplus_{\beta\in\pi}H_{\beta}\longmapsto\bigoplus_{\beta\in\pi}H_{\beta}$.
 
 Conversely, let $H_{\grande}$ be a Hopf algebra with 
 multiplication $\mu_{\grande}$, unit $1$, and antipode $S_{\grande}$, 
 endowed with a group homomorphism
 $\map{\varphi_{\grande}}{\pi}{\Aut(H_{\grande})}{\alpha}{\varphi_{\grande,\alpha}}$.
 Suppose that the following conditions are satisfied.
 \begin{itemize}
  \item There exists a family of sub-coalgebras $\{H_{\alpha}\}_{\alpha\in\pi}$ of  $H_{\grande}$
        such that $H_{\grande}=\bigoplus_{\alpha\in\pi}H_{\alpha}$.
  \item $H_{\alpha}\cdot H_{\beta}\subset H_{\alpha\beta}$ for any $\alpha,\beta\in\pi$.
  \item $1\in H_{1}$.
  \item For any $\alpha,\beta\in\pi$, $\varphi_{\grande,\beta}$ sends $H_{\alpha}\subset H_{\grande}$ to
        $H_{\beta\alpha\beta^{-1}}\subset H_{\grande}$.
  \item For any $\alpha\in\pi$, we have $S_{\grande}(H_{\alpha})=H_{\alpha^{-1}}$
 \end{itemize}
 
\noindent Then, we obtain, in the obvious way, 
a \Talg $H$ such that $H_{\pk}=H_{\grande}$.
  
\section{The outer dual and the inner dual of a \protect\Tcoalg}
 \label{sec:duals}
 We study how to provide a convenient
 notion of a dual for a finite-type \Tcoalg $H$. The easiest way is
 to define a \Talg $\Houter$, the \textit{outer dual of $H$.} However,
 for many purposes, in particular in the construction of the quantum
 double of $H$, it is convenient to introduce a \THcoalg $H^{\innersym}$
 based on the packed form $\Houter_{\pk}$ of $\Houter$, the 
 \textit{inner dual of $H$.} 
 
\paragraph{\scshape The outer dual}
 Let $H$ be a finite-type \Tcoalg. The \textit{outer dual of $H$} 
 is the \Talg $\Houter$ defined as follows.
 For any $\alpha\in\pi$, the \alphath component of $\Houter$ is the 
 dual coalgebra $H^{\ast}_{\alpha}$ of the algebra $H_{\alpha}$. 
 The multiplication of $\Houter$ is given by
 \begin{equation}\label{e:Oratius}
  \bigl\langle \mu_{\alpha,\beta}(f,g),h\bigr\rangle =\langle f\otimes g, \Delta_{\alpha,\beta}(h)\rangle
 \end{equation}
 for any $f\in H^{\ast}_{\alpha}$, $g\in H^{\ast}_{\beta}$ and $h\in H_{\alpha\beta}$, with $\alpha,\beta\in\pi$.
 The unit of $\Houter$ is given by $\varepsilon\in H_{1}^{\ast}\subset\Houter$. 
 The antipode $S^{\outersym}$ of $H^{\outersym}$ is given by
 $S^{\outersym}_{\alpha}\eqdef s^{\ast}_{\alpha^{-1}}$,
 for any $\alpha\in\pi$. Finally, for any $\beta\in\pi$, the  conjugation 
 isomorphism $\psi^{\outersym}_{\beta}$ of $\Houter$ is given by
 $\psi^{\outersym}_{\beta}\eqdef\varphi^{\ast}_{\beta^{-1}}$.
 
\paragraph{\scshape The inner dual}\label{par:inner-dual}
 Since $(\Houter)_{\pk}$ is a Hopf algebra endowed 
 with a group homomorphism
 $\psi_{(\Houter)_{\pk}}\colon\pi\to\Aut\bigl((\Houter)_{\pk}\bigr)$, 
 we can construct the \THcoalg based on $(\Houter)_{\pk}$.
 We call this \THcoalg the \textit{inner dual of $H$} 
 and we denote it $\Hinner$. 
 Explicitly, $\Hinner_{1}=(\Houter)_{\pk}$ is obtained as follows.
 \begin{itemize}
  \item As a coalgebra, $H^{\innersym}_{1}=\bigoplus_{\alpha\in\pi}H^{\ast}_{\alpha}$.     
  \item The multiplication is obtained 
        by~\eqref{e:Oratius}, extending by linearity.
        The unit is given by $\varepsilon^{\innersym}=\varepsilon\in H_{1}^{\ast}\subset \bigoplus_{\alpha\in\pi}H^{\ast}_{\alpha}$.
  \item The antipode is given by 
        $s^{\innersym}_{1}=\sum_{\alpha\in\pi} S^{\outersym}_{\alpha}=\sum_{\alpha\in\pi}s^{\ast}_{\alpha^{-1}}$.
  \item Finally, $\psi_{{\Houter_{\pk},\beta}}=\sum_{\beta\in\pi}\varphi_{\beta^{-1}}^{\ast}$.
\end{itemize}
% 
%\begin{rmk}\label{rmk:conv}
%  If $H$ is a \Tcoalg of finite-type, then
%  $\Conv_{\alpha}(H) =H_{\alpha}\otimes\bigoplus_{\beta\in\pi}H^{\ast}_{\beta}$ for any $\alpha\in\pi$,
%  So, by definition~\eqref{e:conv} of the multiplication 
%  in $\Conv_{\alpha}(H)$, as an algebra 
%  $\Conv_{\alpha}(H)=H_{\alpha}\otimes H^{\innersym}_{1}$.
%\end{rmk}
  
\paragraph{\scshape The coopposite inner dual}\label{p:coopinndoub}
 Given any \Tcoalg $H$, then $\bigl((\Houter)_{\pk}\bigr)^{\cop}$ is 
 the Hopf algebra obtained from $(\Houter)_{\pk}$ by replacing 
 its comultiplication with the new one $\Delta_{\ast}=\Delta^{\innersym,\cop}$ given by
 \begin{equation}\label{e:ref-a-1}
  \bigl\langle \Delta_{\ast}(f),h\otimes k\bigr\rangle = \langle f, kh\rangle
 \end{equation}
 for any $f\in H^{\ast}_{\alpha}\subset\bigoplus_{\beta\in\pi}H_{\beta}^{\ast}$ and $h,k\in H_{\alpha}$, with $\alpha\in\pi$.
 We also need to replace the antipode with the new 
 one given by $s_{\ast}=S^{\outersym,\cop}=(S^{\outersym})^{-1}$. 
 In particular, we have
 $\bigl\langle s_{\ast}(f),h\bigr\rangle = \bigl\langle f, s_{\alpha}^{-1}(h)\bigr\rangle$,
 for any $f\in H_{\alpha}^{\ast}$ and $h\in H_{\alpha^{-1}}$, with $\alpha\in\pi$.
 The \THcoalg based on $\bigl((\Houter)_{\pk}\bigr)^{\cop}$ 
 is called the \textit{coopposite inner dual of $H$} 
 and is denoted $\Hdualcop$. Notice that
 $\varphi_{\Hdualcop,\alpha}=\varphi_{H^{\innersym},\alpha}=\sum_{\beta\in\pi}\varphi^{\ast}_{\beta^{-1}}$, for any $\alpha\in\pi$.

\section{Quasitriangular  \protect\Tcoalgs}\label{fur-Alina}
 The notion of a quasitriangular Hopf algebra~\cite{Drn} is 
 generalized to the case of a \Tcoalg in~\cite{Tur-CPC}.
 A \textit{quasitriangular \Tcoalgs}
 is a \Tcoalg $H$ endowed with a family
 $R=\{R_{\alpha,\beta}=\xi_{(\alpha).i}\otimes\zeta_{(\beta).i}\in H_{\alpha}\otimes H_{\beta}\}_{\alpha,\beta\in\pi}$,
 called \textit{universal \Rmatrix,} such that $R_{\alpha,\beta}$ 
 is invertible for any $\alpha,\beta\in\pi$
 and the following conditions are satisfied.
 \begin{itemize}\begin{subequations}\label{e:R}
   \item For any $\alpha,\beta\in\pi$ and $h\in H_{\alpha\beta}$,
         \begin{equation}\label{e:R-a}     
          R_{\alpha,\beta}\Delta_{\alpha,\beta}(h)=\bigl(\sigma\circ(\varphi_{\alpha^{-1}}\otimes H_{\alpha})\circ
          \Delta_{\alpha\beta\alpha^{-1},\alpha}\bigr)(h)R_{\alpha,\beta}\text{.}
         \end{equation}
   \item For any $\alpha,\beta,\gamma\in\pi$,
         \begin{equation}\label{e:R-b}
          (H_{\alpha}\otimes\Delta_{\beta,\gamma})(R_{\alpha,\beta\gamma})=(R_{\alpha,\gamma})_{1\beta 3}
          (R_{\alpha,\beta})_{12\gamma}\text{,}
         \end{equation}
         where, given two vector spaces $P$ and $Q$ over $\Bbbk$,
         for any $x=p_{i}\otimes q_{i}\in P\otimes Q$ we set
         $x_{1\beta 3}=p_{i}\otimes 1_{\beta}\otimes q_{i}$, and
         $x_{12\gamma}=p_{i}\otimes q_{i}\otimes 1_{\gamma}$.
   \item For any $\alpha,\beta,\gamma\in\pi$,
         \begin{equation}\label{e:R-c}
          (\Delta_{\alpha,\beta}\otimes H_{\gamma})(R_{\alpha\beta,\gamma})=
          \bigl((\varphi_{\beta}\otimes H_{\gamma})(R_{\beta^{-1}\alpha\beta,\gamma})\bigr)_{1\beta 3}
          (R_{\beta,\gamma})_{\alpha 23}\text{,}
         \end{equation}
         where, given two vector spaces $P$ and $Q$,
         for any $x=p_{i}\otimes q_{i}\in P\otimes Q$ we set
         $x_{\alpha 23}= 1_{\alpha}\otimes p_{i}\otimes q_{i}$.
   \item $R$ is \textit{compatible with $\varphi$,} in the sense that,
         for any $\alpha,\beta,\gamma\in\pi$, we have
         \begin{equation}\label{e:R-d}      
          (\varphi_{\alpha}\otimes\varphi_{\alpha})(R_{\beta,\gamma})=R_{\alpha\beta\alpha^{-1},\alpha\gamma\alpha^{-1}}\text{.}
         \end{equation}
  \end{subequations}\end{itemize}
 
\noindent Notice that $(H_{1},R_{1,1})$ is a quasitriangular
Hopf algebra in the usual sense.
 
For any $\alpha, \beta\in\pi$, we introduce the notation
$\tilde{\xi}_{(\alpha).i}\otimes\tilde{\zeta}_{(\beta).i}=\tilde{R}_{\alpha,\beta}=(R^{-1})_{\alpha,\beta}$.
 
\begin{rmk}[Yang-Baxter equation]\label{rmk:YB}
 It is proved in~\cite{Tur-CPC} that, for any $\alpha,\beta,\gamma\in\pi$ we have
  \begin{multline}\label{e:YB}
   (R_{\beta,\gamma})_{\alpha 23}(R_{\alpha,\gamma})_{1\beta 3}(R_{\alpha,\beta})_{12\gamma}\\
   =(R_{\alpha,\beta})_{12\gamma}\bigl((H_{\alpha}\otimes\varphi_{\beta^{-1}})(R_{\alpha,\beta\gamma\beta^{-1}})\bigr)_{1\beta 3}
   (R_{\beta,\gamma})_{\alpha 23}\text{.}
  \end{multline}
  This is an analog of the standard Yang-Baxter equation
  (see, e.g.,~\cite{Kas}).
\end{rmk}
 
\paragraph{\scshape The mirror \protect\Tcoalg}\label{par:mirror}
 Let $H=(H,R)$ be a quasitriangular Hopf algebra (with $R=\xi_{i}\otimes\zeta_{i}$ and 
 $R^{-1}=\tilde{R}=\tilde{\xi}_{i}\otimes\tilde{\zeta}_{i}\/$). By replacing $R$ with 
 $\overline{R}=\sigma(\tilde{R})=\tilde{\zeta}_{i}\otimes\tilde{\xi}_{i}$, 
 we obtain another 
 quasitriangular structure $\overline{H}=(H,\overline{R})$. 
 This means that, in the category of 
 representations of $H$, we replace the braiding $c_{R}$ provided 
 by $R$ by the braiding $c_{R}^{-1}$ provided by $\overline{R}$. 
 When $H$ is a \Tcoalg, the family 
 $\{R^{-1}_{\alpha,\beta}=\tilde{\zeta}_{\alpha.i}\otimes\tilde{\xi}_{\beta.i}\}_{\alpha,\beta\in\pi}$ 
 is not a universal \Rmatrix for $H$. Nevertheless, 
 it is still possible to generalize
 the definition of $\overline{H}$
 in the following way~\cite{Tur-CPC}.
 
Let $H$ be a \Tcoalg. The \Tcoalg $\overline{H}$, called 
\textit{the mirror of $H$}~\cite{Tur-CPC}, is defined as follows.
\begin{itemize}
 \item For any $\alpha\in\pi$, we set $\overline{H}_{\alpha}\eqdef H_{\alpha^{-1}}$.
  \item For any $\alpha,\beta\in\pi$, the component $\overline{\Delta}_{\alpha,\beta}$ of the 
        comultiplication $\overline{\Delta}$
        of $\overline{H}$ is given by
        \begin{equation}\label{e:mirror-temp}
         \overline{\Delta}_{\alpha,\beta}(h)\eqdef \bigl((\varphi_{\beta}\otimes H_{\beta^{-1}})\circ
         \Delta_{\beta^{-1}\alpha\beta,\beta^{-1}}\bigr)(h)\in H_{\alpha^{-1}}\otimes 
         H_{\beta^{-1}}=\overline{H}_{\alpha}\otimes\overline{H}_{\beta}\text{,}
        \end{equation}
        for any $h\in H_{\beta^{-1}\alpha^{-1}}=\overline{H}_{\alpha\beta}$. If we set
        $h'_{\overline{(\alpha)}}\otimes h''_{\overline{(\beta)}}=\overline{\Delta}_{\alpha,\beta}(h)$,
        then~\eqref{e:mirror-temp} can be written in the form
        $h'_{\overline{(\alpha)}}\otimes h''_{\overline{(\beta)}}\eqdef 
         \varphi_{\beta}(h'_{(\beta^{-1}\alpha^{-1}\beta)})\otimes h''_{(\beta^{-1})}$.
        The counit of $\overline{H}$ is given by
        $\varepsilon\in H_{1}^{\ast}=\overline{H}_{1}^{\ast}$.
  \item For any $\alpha\in\pi$, the \alphath component 
        of the antipode $\overline{s}$ of $\overline{H}$ 
        is given by
        $\overline{s}_{\alpha}=\varphi_{\alpha}\circ s_{\alpha^{-1}}$.
  \item Finally, for any $\alpha\in\pi$, we set $\overline{\varphi}_{\alpha}=\varphi_{\alpha}$.
 \end{itemize}
 
 If $H$ is quasitriangular, then $\overline{H}$ 
 is also quasitriangular with universal \Rmatrix $\overline{R}$ 
 given by
 \begin{equation}\label{e:R-mirror}
  \overline{R}_{\alpha,\beta}=\bigl(\sigma(R_{\beta^{-1},\alpha^{-1}})\bigr)^{-1}\in 
  H_{\alpha^{-1}}\otimes H_{\beta^{-1}}=\overline{H}_{\alpha}\otimes\overline{H}_{\beta}
 \end{equation}
 for any $\alpha,\beta\in\pi$. 
 
If, for any $\alpha,\beta\in\pi$, we introduce the notation 
$\overline{\xi}_{(\alpha).i}\otimes\overline{\zeta}_{(\beta).i}\eqdef\overline{R}_{\alpha,\beta}$, 
then we can write~\eqref{e:R-mirror} in the form
$\overline{\xi}_{(\alpha).i}\otimes\overline{\zeta}_{(\beta).i}=
  \tilde{\zeta}_{(\alpha^{-1}).i}\otimes\tilde{\xi}_{(\beta^{-1}).i}$.
 
 Notice that $\overline{\overline{H}}=H$. Notice also that, 
 due to the definition of $\overline{\Delta}$, the mirror of a 
 \THcoalg is not, in general, a \THcoalg. In particular, 
 the mirror of the inner dual of a finite-type \Tcoalg 
 is not a \THcoalg.
 
\section{The quantum double 
         of a finite-type \protect\Tcoalg}\label{sec:QD}
 
\begin{thm}\label{thm:Asterix}
Let $H$ be a finite-type \Tcoalg.
 There exists a unique quasitriangular \Tcoalg $D(H)$,
 the \textup{quantum double of $H$,}
 such that the following conditions are satisfied.
 \begin{itemize}
  \item The \Tcoalgs $\overline{H}$ and $\Hdualcop$ can be 
        embedded into $D(H)$, i.e, there are two morphisms 
        of \Tcoalgs $i\colon\overline{H}\to D(H)$ and
        $j\colon\Hdualcop\to D(H)$
        such that $i_{\alpha}$ and $j_{\alpha}$ are injective for any $\alpha\in\pi$.
  \item For any $\alpha\in\pi$, the linear map
        \begin{equation}\label{e:pro-Asterix}
          p_{\alpha}=\biggl( \Hdualcop_{1} \otimes\overline{H}_{\alpha}
         \xrightarrow{j_{\alpha} \otimes i_{\alpha}}
        D_{\alpha}(H)\otimes D_{\alpha}(H)\xrightarrow{\mu_{\alpha}} D_{\alpha}(H)\biggr)
        \end{equation}
        is bijective \textup{(}where $D_{\alpha}(H)$ is the 
        \alphath component of $D(H)$
        and $\mu_{\alpha}$ is the multiplication in $D_{\alpha}(H)\/$\textup{)}.
  \item For any $\alpha,\beta\in\pi$, the component $R_{\alpha,\beta}$ of the
        \Rmatrix $R$ of $D(H)$ is the image of the
        canonical element of $H_{\alpha^{-1}}\otimes\Hdualcop_{1}$ under the
        embedding $j_{\alpha}\otimes i_{\beta}\colon\Hdualcop_{1}\otimes H_{\alpha^{-1}}\hookrightarrow D_{\alpha}(H)\otimes D_{\beta}(H)$.
 \end{itemize}
\end{thm}

To prove Theorem~\ref{thm:Asterix},
firstly construct a quasitriangular \Tcoalg $D(H)$. 
Then, we will see that $D(H)$ is a solution of the above
universal problem.
 
\begin{rmk}\label{rmk:ref-a-2}
 Notice that, when $\pi=\{1\}$, both $D(H)$ and $D(\overline{H})$ 
 coincide with the standard definition of the quantum double of $H$.
 The reason of the convention that takes the mirror of $H$ in the 
 above description is that we want to be coherent with the
 conventions in~\cite{KasTur}. 
 However, often in the standard case $p_{1}$
 is defined as $\mu_{1}\circ (i_{1}\otimes j_{1})$, reversing the position of
 $i_{1}$ and $j_{1}$. 
 Moreover, some authors identify $D_{1}(H)$ with the vector space 
 $H_{1}^{\ast}\otimes H$.  
 However, since on that point it seems there is no standard 
 convention, here we also follows the notations in~\cite{KasTur},
 apart for the detail that we reversed the order
 of the factors in the tensor product.
\end{rmk}

\paragraph{\scshape Construction of $D(H)$}
 Let $H$ be a finite-type \Tcoalg. 
 We realize $D(H)$ as follows. 
 \begin{itemize}
   \item For any $\alpha\in\pi$, the \alphath component of $D(H)$, 
         denoted $D_{\alpha}(H)$, is, as a vector space,
         \begin{equation*}
          H_{\alpha^{-1}}\otimes \Hdualcop_{\alpha}=H_{\alpha^{-1}}\otimes \Hdualcop_{1}
          = H_{\alpha}\otimes\bigoplus_{\beta\in\pi}H_{\beta}^{\ast}\text{.}
         \end{equation*}
      
         The multiplication in $D_{\alpha}(H)$ is not
         obtained by tensor product of algebras
         of $H_{\alpha^{-1}}$ and $\Hdualcop$.In the 
         sequel, given $h\in H_{\alpha^{-1}}$ and 
         $F\in \Hdualcop_{\alpha}$, the element in $D(H)$
         corresponding  to $h\otimes F$ in $D(H)$ will be denoted $h\circledast F$.
        
         In view of the role played by $\Hdualcop$ in the 
         construction 
         of the quantum double, the 
         Heynemann\hspace{0pt}-\hspace{0pt}Sweedler notation 
         will be reserved for the comultiplication of $\Hdualcop$, 
         not of $H^{\innersym}$, i.e., given 
         $F\in\sum_{\alpha\in\pi}H^{\ast}_{\alpha}$, we set $F'\otimes F''\eqdef\Delta_{\ast}(F)$.
         
         Let $f$ be in $H^{\ast}_{\alpha}$ and let $h$ and $k$ be in $H_{\alpha}$, 
         with $\alpha\in\pi$.
         By $\langle f, h\_ k\rangle$ we denote the linear functional on $H_{\alpha}$ 
         that evaluated at $x\in H_{\alpha}$ gives $\langle f, h x k\rangle$.
         $D_{\alpha}(H)$ is an algebra under the multiplication obtained 
         by setting, for any $h$, $k\in H_{\alpha^{-1}}$, $f\in H^{\ast}_{\gamma}$, 
         and $g\in H^{\ast}_{\delta}$, with $\gamma,\delta\in\pi$,
          \begin{equation}\label{e:D:product}
           (h\circledast f)\,(k\circledast g)\eqdef
           h''_{(\alpha^{-1})}k \circledast f\Bigl\langle g, 
           s^{-1}_{\delta}\bigl(h'''_{(\delta^{-1})}\bigr)\_
           \varphi_{\alpha}\bigl(h'_{(\alpha^{-1}\delta\alpha)}\bigr)\Bigr\rangle\text{.}
         \end{equation}
         The unit of $D_{\alpha}(H)$ is given by $1_{\alpha^{-1}}\circledast\varepsilon$.
         It follows that the canonical embeddings 
         $H_{\alpha^{-1}}, \Hdualcop_{\alpha}\hookrightarrow D_{\alpha}(H)$
         are algebra morphisms and that, for any 
         $h\in H_{\alpha^{-1}}$ and $f\in H^{\ast}_{\gamma}$, we have
         \begin{subequations}\label{e:Golberg-Variationen}
         \begin{gather}
          \label{e:Golberg-Variationen-a}
           (1_{\alpha^{-1}}\circledast f)\,(h\circledast\varepsilon)  = h\circledast f\qquad\text{and}\\
           (h\circledast\varepsilon)\, (1_{\alpha^{-1}}\circledast f) =
           h''_{(\alpha^{-1})}\circledast \bigl\langle f, 
           s_{\gamma}^{-1}(h'''_{(\gamma^{-1})})\_\varphi_{\alpha}(h'_{(\alpha^{-1}\gamma\alpha)})\bigr\rangle\text{.}
          \label{e:Golberg-Variationen-b}
         \end{gather}
         \end{subequations}
%
%         Notice that, for $\pi=\{1\}$, we recover the standard formula of the
%         multiplication of the quantum double of an Hopf algebra, i.e.,
%         \begin{equation*}
%          (h\circledast f)\,(k\circledast g) = h''k\circledast f\bigl\langle g, s^{-1}(h''')\_ h'\bigr\rangle
%         \end{equation*}
%         (for any $h,k\in H_{1}$ and $f,g\in 
%         H^{\ast\scriptscriptstyle\mathrm{cop}}_{1}\/$).
   \item The comultiplication is given by
         \begin{equation}\label{e:D:comultiplication}
          \Delta_{\alpha,\beta}(h\circledast F) =
          \Bigl(\varphi_{\beta}\bigl(h'_{(\beta^{-1}\alpha^{-1}\beta)}\bigr)\circledast F'\Bigr)\otimes 
          \Bigl(h''_{(\beta^{-1})}\circledast F''\Bigr)\text{,}
          \end{equation}
          for any $\alpha,\beta\in\pi$, $h\in \overline{H}_{\alpha\beta}=H_{\beta^{-1}\alpha^{-1}}$, and 
          $F\in \Hdualcop_{\alpha\beta}$ (We recall that $F'\otimes F''=\Delta_{\ast}(F)$, 
          see~\eqref{e:ref-a-1}). 
          The counit is obtained by setting
          $\langle \varepsilon, h\circledast f\rangle = \langle \varepsilon, h\rangle\,\langle f, 1_{\gamma}\rangle$,
          for any $h\in H_{1}$ and $f\in H^{\ast}_{\gamma}$, with $\gamma\in\pi$.
    \item For any $\alpha\in\pi$, the \alphath component of the 
          antipode of $D(H)$ is given by
          \begin{equation}\label{e:D:antipode}
            s_{\alpha}(h\circledast F) = \bigl(\overline{s}_{\alpha}(h)\circledast\varepsilon\bigr)\,
            \bigl(1\circledast s_{\ast}(F)\bigr)
             =\bigl((\varphi_{\alpha}\circ s_{\alpha^{-1}})(h)\circledast\varepsilon\bigr)\,
             \bigl(1_{\alpha}\circledast s_{\ast}(F)\bigr)\text{,}
           \end{equation}
          for any $h\in \overline{H}_{\alpha}=H_{\alpha^{-1}}$ and $F\in\Hdualcop_{\alpha}$, 
          where $s_{\ast}$ is the antipode of $\Hdualcop$ and 
          $\overline{s}_{\alpha}=\varphi_{\alpha}\circ s_{\alpha^{-1}}$ is the antipode of 
          $\overline{H}$.
     \item Finally, for any $\alpha\in\pi$, we set
           \begin{equation}\label{e:D:conj}
            \varphi_{\beta}(h\circledast f)=\varphi_{\beta}(h)\circledast\varphi_{\Hdualcop,\beta}(f)=
            \varphi_{\beta}(h)\circledast\varphi_{\beta^{-1}}^{\ast}(f)\text{,}
           \end{equation}
           for any $h\in H_{\alpha^{-1}}$ and $f\in\Hdualcop_{\gamma}$, with $\gamma\in\pi$.
 \end{itemize}
 
\begin{thm}\label{thm:dual}
  $D(H)$ is a \Tcoalg. Moreover, the multiplication 
  in $D(H)$ is uniquely defined 
  by~\eqref{e:Golberg-Variationen} and the condition
  that the embeddings
  $\overline{H}, \Hdualcop\hookrightarrow D(H)$ are \Tcoalg morphisms.
\end{thm}
 
\begin{proof}
 Firstly, for any $\alpha\in\pi$, we will show that $D_{\alpha}(H)$ 
 is an associative algebra with unit. 
 Then we will show that $\Delta$, defined as above, is multiplicative, 
 i.e., that any $\Delta_{\alpha,\beta}$ is an algebra morphism. 
 After that, we will show that $\varepsilon$ is an algebra morphism. 
 Finally, we will check axioms for the antipode and that 
 the conjugation isomorphisms are compatible with the multiplication.
 
\begin{sentence}{Associativity}
 Let $\alpha$ be in $\pi$. The multiplication defined in~\eqref{e:D:product} 
 is associative if and only if, for any $h, k, l\in H_{\alpha^{-1}}$, 
 $p\in H_{\beta}^{\ast}$, $q\in H_{\gamma}^{\ast}$, and $r\in H_{\delta}^{\ast}$, with $\beta,\gamma,\delta\in\pi$,
   \begin{equation}\label{e:D:product-assoc}
    \bigl((h\circledast p)(k\circledast q)\bigr)(l\circledast r)= (h\circledast p)\bigl((k\circledast q)(l\circledast r)\bigr)
   \end{equation}
 By computing the left-hand side of~\eqref{e:D:product-assoc}, we obtain
   \begin{equation*}\begin{split}
    \ \bigl((h &\circledast p)(k\circledast q)\bigr)(l\circledast r)
    = h'''_{(\alpha^{-1})}k''_{(\alpha^{-1})}l\circledast p\Bigl\langle q, 
    s^{-1}_{\gamma}\bigl(h'''''_{(\gamma^{-1})}\bigr)\_\varphi_{\alpha}\bigl(h'_{(\alpha^{-1}\gamma\alpha)}\bigr)\Bigr\rangle\\
    & \qquad\qquad\qquad \Bigl\langle r, s^{-1}_{\delta}\bigl(h''''_{(\delta^{-1})}k'''_{(\delta^{-1})}\bigr)\_
    \varphi_{\alpha}\bigl(h''_{(\alpha^{-1}\delta\alpha)}k'_{(\alpha^{-1}\delta\alpha)}\bigr)\Bigr\rangle\\
          \intertext{(by the antimultiplicativity 
          of $s$ and the multiplicativity of $\varphi\/$)}
    & = h'''_{(\alpha^{-1})}k''_{(\alpha^{-1})}l\circledast p\Bigl\langle q, 
    s^{-1}_{\gamma}\bigl(h'''''_{(\gamma^{-1})}\bigr)\_\varphi_{\alpha}\bigl(h'_{(\alpha^{-1}\gamma\alpha)}\bigr)\Bigr\rangle\\
    &\qquad\qquad\qquad\Bigl\langle r, s^{-1}_{\delta}\bigl(k'''_{(\delta^{-1})}\bigr)
    s^{-1}_{\delta}\bigl(h''''_{(\delta^{-1})}\bigr)\_\varphi_{\alpha}\bigl(h''_{(\alpha^{-1}\delta\alpha)}\bigr)
    \varphi_{\alpha}\bigl(k'_{(\alpha^{-1}\delta\alpha)}\bigr)\Bigr\rangle\text{,}    
   \end{split}\end{equation*}
   while, by computing the right-hand side, we obtain
   \begin{equation*}\begin{split}
    \quad\ (h \circledast p&)\bigl(\bigl(k\circledast q)(l\circledast r)\bigr)=
     h''_{(\alpha^{-1})}k''_{(\alpha^{-1})}l\circledast p\biggl\langle q\Bigl\langle r, 
    s^{-1}_{\delta}\bigl(k'''_{(\delta^{-1})}\bigr)\_\varphi_{\alpha}\bigl(k'_{(\alpha^{-1}\delta\alpha)}\bigr)\Bigr\rangle,\\
    &\qquad\qquad\qquad\qquad\qquad\qquad
    s^{-1}_{\gamma\delta}\bigl(h'''_{(\delta^{-1}\gamma^{-1})}\bigr)\_\varphi_{\alpha}\bigl(h'_{(\alpha^{-1}\gamma\delta\alpha)}\bigr)\biggr\rangle
           \intertext{(by the anticomultiplicativity of 
           $s$ and the comultiplicativity of $\varphi\/$)}
    & = h'''_{(\alpha^{-1})}k''_{(\alpha^{-1})}l\circledast p\Bigl\langle q, 
    s^{-1}_{\gamma}\bigl(h'''''_{(\gamma^{-1})}\bigr)\_\varphi_{\alpha}\bigl(h'_{(\alpha^{-1}\gamma\alpha)}\bigr)\Bigr\rangle\\
    & \qquad\qquad\qquad
    \Bigl\langle r, s^{-1}_{\delta}\bigl(k'''_{(\delta^{-1})}\bigr)s^{-1}_{\delta}\bigl(h''''_{(\delta^{-1})}\bigr)\_
    \varphi_{\alpha}\bigl(h''_{(\alpha^{-1}\delta\alpha)}\bigr)\varphi_{\alpha}\bigl(k'_{(\alpha^{-1}\delta\alpha)}\bigr)\Bigr\rangle\text{.}
   \end{split}\end{equation*}  
 \end{sentence}
 
  \begin{sentence}{Unit}
   Let $\alpha$ be in $\pi$. For any $h\in H_{\alpha^{-1}}$ 
   and $f\in H^{\ast}_{\gamma}$, with $\gamma\in\pi$, we have
   \begin{equation*}
    (1_{\alpha^{-1}}\circledast\varepsilon)\,(h\circledast f)=1_{\alpha^{-1}}h\circledast \varepsilon f=h\circledast f
   \end{equation*}
   and
   \begin{equation*}
    (h\circledast f)\,(1_{\alpha}\circledast\varepsilon)= h''_{(\alpha^{-1})}1_{\alpha^{-1}}\circledast f\varepsilon\bigl\langle\varepsilon, 
    s^{-1}_{1}(h''')\bigr\rangle\bigl\langle\varepsilon, \varphi_{\alpha}(h')\bigr\rangle=h\circledast f\text{,}
   \end{equation*}
   where the fact that both $s_{1}$ and $\varphi_{\alpha}$ commute with $\varepsilon$.
  \end{sentence}
    
  \begin{sentence}{Multiplicativity of $\Delta$}
   Let us prove that $\Delta_{\alpha,\beta}$ 
   is an algebra morphism for any $\alpha,\beta\in\pi$.
   Since $\Delta_{\alpha,\beta}$ obviously preserves the unit, we only need to 
   prove that, for any 
   $h, k\in H_{\beta^{-1}\alpha^{-1}}$, $f\in H^{\ast}_{\gamma}$ and $g\in H^{\ast}_{\delta}$, with $\gamma,\delta\in\pi$, 
   we have
   \begin{equation}\label{e:D:temp-ass-delta}
    \Delta_{\alpha,\beta}\bigl((h\circledast f)(k\circledast g)\bigr)=
    \Delta_{\alpha,\beta}(h\circledast f)\Delta_{\alpha,\beta}(k\circledast g)\text{.}
   \end{equation}
   This is proved by evaluating both terms in~\eqref{e:D:temp-ass-delta}
   against the general term $p\otimes x\otimes q\otimes y$ ($p\in H^{\ast}_{\alpha^{-1}}$, 
   $q\in H^{\ast}_{\beta^{-1}}$, and $x,y\in H_{\gamma\delta}$).
  \end{sentence}
   
  \begin{sentence}{Multiplicativity of $\varepsilon$}
   For any $h, k\in H_{1}$, 
   $f\in H^{\ast}_{\gamma}$, and $g\in H^{\ast}_{\delta}$, with $\gamma,\delta\in\pi$, we have
   \begin{equation*}
    \langle\varepsilon, h\circledast f\rangle\,\langle\varepsilon, k\circledast f\rangle = \langle 
    \varepsilon, h\rangle\,\langle f, 1_{\gamma}\rangle\,\langle\varepsilon, k\rangle\,\langle f, 1_{\delta}\rangle
   \end{equation*}
   and
   \begin{equation*}\begin{split}
   \bigl\langle\varepsilon,  (h\circledast f)(k\circledast g)\bigr\rangle & =
   \Bigl\langle\varepsilon, h'' k \circledast f\bigl\langle g, s^{-1}_{\delta}(h'''_{(\delta^{-1})})\_ h'_{(\delta)}\bigr\rangle\Bigr\rangle\\
   & = \langle\varepsilon, h''\rangle\,\langle\varepsilon, k\rangle\,\bigl\langle f, 1_{\gamma}\bigr\rangle\,
   \bigl\langle g, s^{-1}_{(\delta)}(h'''_{(\delta^{-1})})h'_{(\delta)}\bigr\rangle\\
   & = \langle\varepsilon, k\rangle\,\bigl\langle f, 1_{\gamma}\bigr\rangle\,\bigl\langle g, 
           s^{-1}_{\delta}(h''_{(\delta^{-1})})h'_{(\delta)}\bigr\rangle \\
   & =  \langle\varepsilon, k\rangle\,\bigl\langle f, 1_{\gamma}\bigr\rangle\,
    \bigl\langle g, \langle\varepsilon, h\rangle 1_{\delta}\bigr\rangle=\langle\varepsilon,h\rangle\,\langle f,1_{\gamma}\rangle\,\langle\varepsilon, k\rangle\,\langle g, 1_{\delta}\rangle\text{.}
   \end{split}\end{equation*}
   This proves that $\varepsilon$ is multiplicative. 
   Moreover, since $\varepsilon$ is obviously unitary, 
   it is an algebra homomorphism.
  \end{sentence}
    
  \begin{sentence}{Antipode}
   Let $h$ be in $H_{1}$ and let $f$ 
   be in $H^{\ast}_{\gamma}$, with $\gamma\in\pi$. We have
   \begin{equation*}\begin{split}
    \bigl(h\circledast & f\bigr)'_{(\alpha)}s_{\alpha^{-1}}\bigl((h\circledast f)''_{(\alpha^{-1})}\bigr)\\
    & = \bigl(\varphi_{\alpha^{-1}}(h'_{(\alpha^{-1})})\circledast f'\bigr)
    \bigl((\varphi_{\alpha^{-1}}\circ s_{\alpha})(h''_{(\alpha)})\circledast\varepsilon\bigr)
    \bigl(1_{\alpha^{-1}}\circledast s_{\ast}(f'')\bigr)\\
    & = \Bigl(\varphi_{\alpha^{-1}}\bigl(h'_{(\alpha^{-1})}s_{\alpha}(h''_{(\alpha)})\bigr)\circledast  f'\Bigr)
    \bigl(1_{\alpha^{-1}}\circledast s_{\ast}(f'')\bigr)\\
    & = \langle\varepsilon, h\rangle\bigl(1_{\alpha^{-1}}\circledast f'\bigr)\bigl(1\circledast s_{\ast}(f'')\bigr)
      = \langle\varepsilon, h\rangle 1\circledast f' s_{\ast}(f'')\\
    & = \langle\varepsilon, h\rangle\langle f, 1_{\gamma}\rangle 1_{\alpha^{-1}}\circledast\varepsilon = \langle\varepsilon, h\circledast f\rangle 1_{\alpha^{-1}}\circledast\varepsilon
   \end{split}\end{equation*}
   and
   \begin{equation*}\begin{split}
    s_{\alpha^{-1}}& \bigl((h\circledast f)'_{(\alpha^{-1})}\bigr)\bigl(h\circledast f \bigr)''_{(\alpha)}
    = s_{\alpha^{-1}}\bigl(\varphi_{\alpha}(h'_{(\alpha)})\circledast f'\bigr)\bigl(h''_{(\alpha^{-1})}\circledast
    f''\bigr)\\ 
    & =\bigl(s_{\alpha}(h'_{\alpha})\circledast\varepsilon \bigr)\bigl(1_{\alpha}\circledast s_{\ast}(f')\bigr)\bigl(h''_{(\alpha^{-1})}
    \circledast f''\bigr)\\
    & = \langle f, 1_{\gamma}\rangle \bigl(s_{\alpha}(h'_{(\alpha)})\circledast\varepsilon\bigr)
    \bigl(h''_{(\alpha^{-1})}\circledast\varepsilon\bigr)= \langle f, 1_{\gamma}\rangle s_{\alpha}(h''_{(\alpha)})h''''_{(\alpha^{-1})}\circledast 
    \langle\varepsilon, h'\rangle\langle\varepsilon, h''\rangle\varepsilon\\ & =
    \langle\varepsilon, h\rangle\langle f, 1_{\gamma}\rangle 1_{\alpha}\circledast\varepsilon = 
    \langle\varepsilon, h\circledast f\rangle 1_{\alpha}\circledast\varepsilon\text{.}   
   \end{split}\end{equation*}
  \end{sentence}
    
  \begin{sentence}{Conjugation}
   Let us check that $\varphi_{\beta}^{\alpha}$
   is an algebra isomorphism for any $\alpha,\beta\in\pi$. Since $\varphi_{\beta}^{\alpha}$ is
   obviously bijective and preserve the unit,
   we only need to show that it is compatible with the multiplication, 
   i.e., that, for any $h, k\in H_{\alpha^{-1}}$, 
   $f\in H_{\gamma}^{\ast}$, and $g\in H_{\delta}^{\ast}$, with $\gamma,\delta\in\pi$, we have
   \begin{equation}\label{e:because-music}
    \varphi_{\beta}(h\circledast f)\varphi_{\beta}(k\circledast g)=
    \varphi_{\beta}\bigl((h\circledast f)(k\circledast g)\bigr)\text{.}
   \end{equation}
   This is proved by evaluating both sides in~\eqref{e:because-music} 
   against the general term $p\otimes x$ 
   ($x\in H_{\beta\gamma\delta\beta^{-1}}$ and $p\in H^{\ast}_{\beta\alpha^{-1}\beta^{-1}}$). 
   
   This concludes the proof that $D(H)$ is a \Tcoalg.
  \end{sentence}
  
  \begin{sentence}{Embeddings}
   It only remains to show that~\eqref{e:Golberg-Variationen}, 
   together with the request that the canonical embeddings
   $\overline{H}, \Hdualcop\hookrightarrow D(H)$ are \Tcoalg morphisms, 
   uniquely determinates the multiplication in $D(H)$. 
   Let $\alpha$ be in $\pi$.
   For any $h,k\in H_{\alpha^{-1}}$, $f\in H^{\ast}_{\gamma}$ and $g\in H^{\ast}_{\delta}$ we have
   \begin{equation*}\begin{split}
    \ (h\circledast& f)\,(k\circledast g) =\\
    \intertext{(by~\eqref{e:Golberg-Variationen-a})}
    & = (1_{\alpha^{-1}}\circledast f)\,(h\circledast\varepsilon)\,(1_{\alpha^{-1}}\circledast g)\,(k\circledast\varepsilon)\\
    \intertext{(by~\eqref{e:Golberg-Variationen-b})} 
    & = (1_{\alpha^{-1}}\circledast f)\,\Bigl(h''_{(\alpha^{-1})}\circledast\bigl\langle g, s^{-1}_{\delta}
      (h'''_{(\delta^{-1})})\_\varphi_{\alpha}(h'_{(\alpha^{-1}\delta\alpha)})\bigr\rangle\Bigr)\,(k\circledast\varepsilon)\\
    \intertext{(by~\eqref{e:Golberg-Variationen-a})}
    & = (1_{\alpha^{-1}}\circledast f)\,\Bigl(1_{\alpha^{-1}}\circledast\bigl\langle g, s^{-1}_{\delta}
      (h'''_{(\delta^{-1})})\_\varphi_{\alpha}(h'_{(\alpha^{-1}\delta\alpha)})\bigr\rangle\Bigr)\,
      (h''_{(\alpha^{-1})}\circledast\varepsilon)\,(k\circledast\varepsilon)\\
    \intertext{(since the canonical embedding of 
    both $H_{\alpha^{-1}}$ and $\Hdualcop_{\alpha}$ in $D_{\alpha}(H)$ are algebra morphisms)}
    & = \Bigl(1_{\alpha^{-1}}\circledast f\bigl\langle g, s^{-1}_{\delta}
      (h'''_{(\delta^{-1})})\_\varphi_{\alpha}(h'_{(\alpha^{-1}\delta\alpha)})\bigr\rangle\Bigr)\,(h''_{(\alpha^{-1})}k\circledast\varepsilon)
    \intertext{(again by~\eqref{e:Golberg-Variationen-a})}
    & = h''_{(\alpha^{-1})}k \circledast f\Bigl\langle g, s^{-1}_{\delta}\bigl(h'''_{(\delta^{-1})}\bigr)\_
      \varphi_{\alpha}\bigl(h'_{(\alpha^{-1}\delta\alpha)}\bigr)\Bigr\rangle\text{.}
   \end{split}\end{equation*}
   This concludes the proof of the theorem.
  \end{sentence}
  
   Coassociativity, compatibility 
 between the comultiplication and the counit,
 compatibility between the comultiplication and the conjugation, 
 and, finally, the fact that $\varphi$ is a group homomorphism follow by
 observing that $D(H)$ has the same comultiplication, counit, and 
 The omitted computations are the same needed 
 conjugation of the tensor product of \Tcoalgs $\overline{H}\otimes\Hdualcop$.
 \end{proof}
 
\paragraph{\scshape Quasitriangular structure of the quantum double}
 To prove that $D(H)$ is quasitriangular, 
 we need some preliminary results.
 For any $\alpha\in\pi$, we set $n_{\alpha}=\dim H_{\alpha}$. 
 Let $(e_{\alpha.i})_{i=1,\ldots,n_{\alpha}}$ be a basis of $H_{\alpha}$ as 
 a vector space and let $(e^{\alpha.i})_{i=1,\ldots,n_{\alpha}}$ 
 be the dual basis of $(e_{\alpha.i})_{i=1,\ldots,n_{\alpha}}$. We set
   \begin{equation}\label{e:D:R}
   R_{\alpha,\beta}\eqdef e_{\alpha^{-1}.i}\circledast\varepsilon\otimes 1_{\beta^{-1}}\circledast e^{\alpha^{-1}.i}
   \in D_{\alpha}(H)\otimes D_{\beta}(H)
  \end{equation}
  and $\tilde{R}_{\alpha,\beta}\eqdef s_{\alpha}(e_{\alpha.i})\circledast\varepsilon\otimes 
   1_{\beta^{-1}}\circledast e^{\alpha.i}$.
 
 \begin{lemma}\label{D:R}
  For any $\alpha,\beta\in\pi$, both $R_{\alpha,\beta}$ and $\tilde{R}_{\alpha,\beta}$ 
  are independent of the choice of the bases. 
  Moreover, $\tilde{R}_{\alpha,\beta}$ is the inverse of $R_{\alpha,\beta}$ 
  in the algebra $D_{\alpha}(H)\otimes D_{\beta}(H)$.
 \end{lemma}
 
 \begin{proof}
  Let $\Cv_{\alpha,\beta}(H)$ be the subspace of $D_{\alpha}(H)\otimes D_{\beta}(H)$ 
  generated by the elements
  $h\circledast \varepsilon \otimes 1_{\beta^{-1}}\circledast F$ with $h\in H_{\alpha^{-1}}$ and $F\in\Hdualcop_{\alpha}$.
  Clearly, $\Cv_{\alpha,\beta}(H)$ is a subalgebra of $D_{\alpha}(H)\otimes D_{\beta}(H)$.
  Moreover, 
  we have an equivalence of algebras 
  $\Conv_{\alpha^{-1}}(H)\to\Cv_{\alpha,\beta}(H)\colon h\otimes F\mapsto (h\circledast \varepsilon)\otimes (1_{\beta^{-1}} \circledast F)$, 
  for any $h\in H_{\alpha^{-1}}$ and $F\in\Hdualcop_{\alpha}$. 
  In particular, $R_{\alpha,\beta}$ is the image of the identity morphism 
  $e_{\alpha^{-1}.i}\otimes e^{\alpha^{-1}.i}$ of 
  $H_{\alpha^{-1}}$ under this isomorphism. Since the 
  \alphath component of the antipode $s$ of $H$ is the inverse, 
  in the algebra $\Conv_{\alpha^{-1}}(H)$, of $e_{\alpha^{-1}.i}\otimes e^{\alpha^{-1}.i}$, 
  also $R_{\alpha,\beta}$ is invertible. Since $s_{\alpha}$ can be 
  represented as $s(e_{\alpha.i})\otimes e^{\alpha.i}$ and $\tilde{R}_{\alpha,\beta}$ 
  is the image of $s_{\alpha}$ under the above isomorphism, we conclude 
  that $\tilde{R}_{\alpha,\beta}$ is the inverse of $R_{\alpha,\beta}$.
 \end{proof}
     
 \begin{rmk}
  Let $h$ be in $H_{\alpha}$, with $\alpha\in\pi$, and suppose
  $\alpha_{1}\alpha_{2}\cdots\alpha_{n}=\alpha$ for certain $\alpha_{1},\alpha_{2}\ldots,\alpha_{n}\in\pi$.
  By observing that $h=h^{i}e_{\alpha.i}=\langle e^{\alpha.i}, h\rangle e_{\alpha.i}$ and by 
  the linearity of $\Delta$ we have
  \begin{equation}\label{l:basis}
   h'_{(\alpha_{1})}\otimes h''_{(\alpha_{2})}\cdots\otimes {h'}{}^{n}_{(\alpha_{n})} 
   =\langle e^{\alpha.i}, h\rangle (e_{\alpha.i})'_{(\alpha_{1})}\otimes (e_{\alpha.i})''_{(\alpha_{2})}\otimes\cdots\otimes 
   (e_{\alpha.i})'{}^{n}_{(\alpha_{n})}\text{.}
  \end{equation}
 \end{rmk}
 
 \begin{thm}\label{thm:DH}
  $D(H)$ is quasitriangular with universal \Rmatrix
  given by~\eqref{e:D:R} for any $\alpha,\beta\in\pi$.
 \end{thm}

\begin{proof}
 $R_{\alpha,\beta}$ is well defined and invertible by Lemma~\ref{D:R}. 
 We still need to check the four relations~\eqref{e:R}.

 \begin{sentence}{Relation~\eqref{e:R-a}}
  Let $\alpha,\beta$, and $\gamma$ be $\pi$. 
  Given $h\in H_{\beta^{-1}\alpha^{-1}}$ and $f\in H^{\ast}_{\gamma}$, we have
%  $ R_{\alpha,\beta} \Delta_{\alpha,\beta}(h\circledast f) = \bigl(e_{\alpha^{-1}.i}\circledast\varepsilon\otimes 1_{\beta^{-1}}\circledast e^{\alpha^{-1}.i}
%   \bigr) \bigl(\varphi_{\beta}(h'_{(\beta^{-1}\alpha^{-1}\beta)})\circledast f'\otimes h''_{(\beta^{-1})}\circledast f''\bigr)
%    =(e_{\alpha^{-1}.i})''_{(\alpha^{-1})}\varphi_{\beta}(h'_{(\beta^{-1}\alpha\beta)})\circledast\Bigl\langle f',
%   s^{-1}_{\gamma}\bigl((e_{\alpha^{-1}.i})'''_{(\gamma^{-1})}\bigr)\_\varphi_{\alpha}
%   \bigl((e_{\alpha^{-1}.i})'_{(\alpha^{-1}\gamma\alpha)}\bigr)\Bigr\rangle\otimes h''_{(\beta^{-1})}\circledast 
%   e^{\alpha^{-1}.i}f''$. 
  \begin{equation*}\begin{split} 
  R_{\alpha,\beta} \Delta_{\alpha,\beta} &(h\circledast f)  = \bigl((e_{\alpha^{-1}.i}\circledast\varepsilon)\otimes( 1_{\beta^{-1}}\circledast e^{\alpha^{-1}.i})
   \bigr) \bigl((\varphi_{\beta}(h'_{(\beta^{-1}\alpha^{-1}\beta)})\circledast f')\otimes( h''_{(\beta^{-1})}\circledast f''\bigr))\\
   & =(e_{\alpha^{-1}.i})''_{(\alpha^{-1})}\varphi_{\beta}(h'_{(\beta^{-1}\alpha\beta)})\circledast\\
   &\qquad \circledast\Bigl\langle f',
   s^{-1}_{\gamma}\bigl((e_{\alpha^{-1}.i})'''_{(\gamma^{-1})}\bigr)\_\varphi_{\alpha}
   \bigl((e_{\alpha^{-1}.i})'_{(\alpha^{-1}\gamma\alpha)}\bigr)\Bigr\rangle\otimes h''_{(\beta^{-1})}\circledast 
   e^{\alpha^{-1}.i}f''\text{.}
  \end{split}\end{equation*}    
   Since the action of 
  $\Bigl(\sigma\circ \bigl(\varphi_{\alpha}\otimes D_{\alpha}(H)\bigr)\circ \Delta_{\alpha\beta\alpha^{-1},\alpha}\Bigr)(\_)$
  on $h\circledast f$ is given by
  \begin{equation*}\begin{split}
   h\otimes f& \overmapsto[3.5pc]{\Delta_{\alpha\beta\alpha^{-1},\alpha}}\varphi_{\alpha}(h'_{(\beta^{-1})})\circledast 
   f'\otimes h''_{(\alpha^{-1})}\circledast f''\overmapsto[4.5pc]{\varphi_{\alpha^{-1}}\otimes D_{\alpha}(H)}\\ 
   & h'_{(\beta^{-1})}\circledast\varphi_{\alpha}^{\ast}(f')\otimes h''_{(\alpha^{-1})}\circledast f''\overmapsto[2pc]{\sigma} 
   h''_{(\alpha^{-1})}\circledast f''\otimes h'_{(\beta^{-1})}\circledast\varphi_{\alpha}^{\ast}(f')\text{,}
  \end{split}\end{equation*}
  we have
\begin{equation*}\begin{split}
  \biggl(\Bigl(\sigma \circ & \bigl(\varphi_{\alpha}\otimes D_{\alpha}(H)\bigr)\circ 
  \Delta_{\alpha\beta\alpha^{-1},\alpha}\Bigr)(h\otimes f)\biggr)R_{\alpha,\beta}\\
  & = \bigl(h''_{(\alpha^{-1})}\circledast f''\otimes h'_{(\beta^{-1})}\circledast\varphi_{\alpha}^{\ast}(f')\bigr)
  \bigl(e_{\alpha^{-1}.i}\circledast\varepsilon\otimes 1_{\beta^{-1}}\circledast e^{\alpha^{-1}.i}\bigr)\\
  & = h''''_{(\alpha^{-1})}e_{\alpha^{-1}.i}\circledast f''\otimes h''_{(\beta^{-1})}\circledast\varphi_{\alpha}^{\ast}(f')
  \bigl\langle e^{\alpha^{-1}.i}, 
  s_{\alpha^{-1}}^{-1}(h'''_{(\alpha)})\_\varphi_{\beta}(h'_{(\beta^{-1}\alpha^{-1}\beta)})\bigr\rangle\text{.}
 \end{split}\end{equation*}
  
  Relation~\eqref{e:R-a} is proved by observing that
  evaluating the two expressions above against the tensor  
  $H_{\alpha^{-1}}\otimes H^{\ast}_{\gamma}\otimes H_{\beta^{-1}}\otimes\langle\cdot,x\rangle$ 
  (for a generic $x\in H_{\alpha^{-1}\gamma}$\/)
  we found the same result.
 \end{sentence}
 
\begin{sentence}{Relation~\eqref{e:R-b}}
  For any $\alpha,\beta,\gamma\in\pi$, we have
  \begin{equation*}\begin{cases}
   (R_{\alpha,\gamma})_{1\beta 3}= (e_{\alpha^{-1}.i}\circledast\varepsilon)\otimes(1_{\beta^{-1}}\circledast\varepsilon)\otimes (1_{\gamma^{-1}}\circledast e^{\alpha^{-1}.i}) & \ \\
   (R_{\alpha,\beta})_{12\gamma}= (e_{\alpha^{-1}.i}\circledast\varepsilon)\otimes(1_{\beta^{-1}}\circledast e^{\alpha^{-1}.i})\otimes (1_{\gamma^{-1}}\circledast\varepsilon) & \ 
  \end{cases}\end{equation*}
  and so we have
   $(R_{\alpha,\gamma})_{1\beta 3}(R_{\alpha,\beta})_{12\gamma} = e_{\alpha^{-1}.i}e_{\alpha^{-1}.j}\circledast\varepsilon\otimes
   1_{\beta^{-1}}\circledast e^{\alpha^{-1}.j}\otimes
   1_{\gamma^{-1}}\otimes e^{\alpha^{-1}.i}$.
  Since 
  $R_{\alpha,\beta\gamma}=e_{\alpha^{-1},i}\circledast\varepsilon\otimes 1_{(\beta\gamma)^{-1}}\circledast e^{\alpha^{-1}.i}$, we have
  $\bigl(D_{\alpha}(H)\otimes\Delta_{\beta,\gamma}\bigr)(R_{\alpha,\beta\gamma})= e_{\alpha^{-1}.i}\circledast\varepsilon \otimes 
   1_{\beta^{-1}}\circledast (e^{\alpha^{-1}.i})'\otimes 1_{\gamma^{-1}}\circledast (e^{\alpha^{-1}.i})''$.
  So, we only need to prove
  \begin{equation}\label{e:cup-of-tea}
   e_{\alpha^{-1}.i}e_{\alpha^{-1}.j}\otimes e^{\alpha^{-1}.j}\otimes e^{\alpha^{-1}.i}=
   e_{\alpha^{-1}.i}\otimes (e^{\alpha^{-1}.i})'\otimes (e^{\alpha^{-1}.i})''\text{.}
  \end{equation}
  If we evaluate both sides 
  of~\eqref{e:cup-of-tea} against the tensor 
  $f\otimes H^{\ast}_{\alpha^{-1}}\otimes H^{\ast}_{\alpha^{-1}}$ (for a generic $f\in H^{\ast}_{\alpha}$\/),
  then, by~\eqref{l:basis}, we found the same result.
 \end{sentence}

 \begin{sentence}{Relation~\eqref{e:R-c}}
 Let $\alpha$, $\beta$, and $\gamma$ be in $\pi$. 
 Observing that
 $R_{\beta^{-1}\alpha\beta,\gamma}=e_{\beta^{-1}\alpha^{-1}\beta.i}\circledast\varepsilon\otimes 1_{\gamma^{-1}}\circledast 
  e^{\beta^{-1}\alpha^{-1}\beta.i}$,
 we obtain
 \begin{equation*} \Bigl(\bigl(\varphi_{\beta}\otimes D_{\gamma}(H)\bigr)(R_{\beta^{-1}\alpha\beta,\gamma})\Bigr)_{1\beta 3}
  = \varphi_{\beta}(e_{\beta^{-1}\alpha^{-1}\beta.i})\circledast\varepsilon\otimes 1_{\beta^{-1}}\circledast
  \varepsilon\otimes 1_{\gamma^{-1}}\circledast e^{\beta^{-1}\alpha^{-1}\beta.i}\text{.}\end{equation*}
 Observing that
 $(R_{\beta,\gamma})_{\alpha 23}= 1_{\alpha^{-1}}\circledast\varepsilon\otimes e_{\beta^{-1}.i}\circledast\varepsilon\otimes 
  1_{\gamma^{-1}}\circledast e^{\beta^{-1}.i}$,
 we obtain
 $\Bigl(\bigl(\varphi_{\beta} \otimes D(H)_{\gamma}\bigr)(R_{\beta^{-1}\alpha\beta,\gamma})\Bigr)_{1\beta 3}(R_{\beta,\gamma})_{\alpha 23} = \varphi_{\beta}(e_{\beta^{-1}\alpha^{-1}\beta.i})\circledast\varepsilon\otimes e_{\beta^{-1}.j}\circledast\varepsilon\otimes 
  1_{\gamma^{-1}}\circledast e^{\beta^{-1}\alpha^{-1}\beta.i}e^{\beta^{-1}.j}$.
 Finally, observing that $R_{\alpha\beta,\gamma}=e_{(\alpha\beta)^{-1},i}\circledast\varepsilon\otimes 1_{\gamma^{-1}}\circledast 
 e^{(\alpha\beta)^{-1}.i}$, we obtain
 $\bigl(\Delta_{\alpha,\beta}\otimes D_{\gamma}(H)\bigr)(R_{\alpha\beta,\gamma})
  =\varphi_{\beta}\bigl((e_{(\alpha\beta)^{-1}.i})'_{(\beta^{-1}\alpha^{-1}\beta)}\bigr)\circledast\varepsilon\otimes     
  (e_{(\alpha\beta)^{-1}.i})''_{(\beta^{-1})}\circledast\varepsilon\otimes 1_{\gamma^{-1}}\circledast 
  e^{(\alpha\beta)^{-1}.i}$.
 So, to prove~\eqref{e:R-c}, we only need to show the equality
 \begin{equation}\label{e:ravi}\begin{split}
  e_{\beta^{-1}\alpha^{-1}\beta.i}\otimes& e_{\beta^{-1}.j}\otimes e^{\beta^{-1}\alpha^{-1}\beta.i}e^{\beta^{-1}.j}\\
  & =(e_{(\alpha\beta)^{-1}.i})'_{(\beta^{-1}\alpha\beta)}\otimes (e_{(\alpha\beta)^{-1}.i})''_{(\beta^{-1})}\otimes 
  e^{(\alpha\beta)^{-1}.i}\text{.}
 \end{split}\end{equation}
 If we evaluate both sides of~\eqref{e:ravi} against the tensor 
 $H_{\beta^{-1}\alpha^{-1}\beta}\otimes H_{\beta^{-1}}\otimes\langle\cdot,x\rangle$ 
 (for a generic $x\in H_{\beta^{-1}\alpha^{-1}}$\/), then
 by~\eqref{l:basis} we found the same result.
 \end{sentence}
 
 \begin{sentence}{Relation~\eqref{e:R-d}}
  Given $\alpha, \beta,\gamma\in\pi$, we have
  $(\varphi_{\gamma}\otimes\varphi_{\gamma})(R_{\alpha,\beta})
   =(\varphi_{\gamma}\otimes\varphi_{\gamma})(e_{\alpha^{-1}.i}\circledast\varepsilon\otimes 1_{\beta^{-1}}\circledast e^{\alpha^{-1}.i})
   = \varphi_{\gamma}(e_{\alpha.i})\circledast\varepsilon\otimes 1_{\beta^{-1}}\circledast\varphi_{\gamma}(e^{\alpha^{-1}.i})$.
  Now,  $\varphi_{\gamma}$ is a linear isomorphism, so 
  $\bigl(\varphi_{\gamma}(e_{\alpha^{-1}.i})\bigr)_{i=1,\ldots,n_{\alpha}}$ 
  is a basis of $H_{\gamma\alpha^{-1}\gamma^{-1}}$, and
  $\bigl(\varphi^{\ast}_{\gamma^{-1}}(e^{\alpha.i})\bigr)_{i=1,\ldots,n_{\alpha}}$ is its dual basis.
  So, by Lemma~\ref{D:R}, $R_{\gamma\alpha\gamma^{-1},\gamma\beta\gamma^{-1}}= (\varphi_{\gamma}\otimes\varphi_{\gamma})(R_{\alpha,\beta})$.
 \end{sentence}

 This concludes the proof of the theorem.
\end{proof}

\paragraph{\scshape Proof of Theorem~\ref{thm:Asterix}} 
We need a preliminary lemma.

\begin{lemma}\label{l:Asterix}
 Let $T$ be any quasitriangular \Tcoalg with \Rmatrix $R$.
 For any $\alpha,\beta\in\pi$ and any $x\in T_{\alpha}$ we have
 \begin{equation}\label{e:Asterix}
  s^{-1}_{\beta}(x'''_{(\beta^{-1})})\xi_{(\beta).i}x'_{(\beta)}\otimes\zeta_{(\beta^{-1}\alpha\beta)}x''_{(\beta^{-1}\alpha\beta)}
  =\xi_{(\beta).i}\otimes\varphi_{\beta^{-1}}(x)\zeta_{(\beta^{-1}\alpha\beta)}\text{.}
 \end{equation}
\end{lemma}

\begin{proof}
 By~\eqref{e:R-a}, we have
  $\xi_{(\beta).i}x'_{(\beta)}\otimes\zeta_{(\beta^{-1}\alpha\beta).i}x''_{(\beta^{-1}\alpha\beta)}\otimes x'''_{(\beta^{-1})}=
  x''_{(\beta^{-1})}\xi_{(\beta).i}\otimes\varphi_{\beta^{-1}}(x'_{(\alpha)})\zeta_{(\beta^{-1}\alpha\beta).i}\otimes x'''_{(\beta^{-1})}$.
 If, on both sides, we apply $T_{\beta}\otimes T_{\beta^{-1}\alpha\beta}\otimes s^{-1}_{\beta}$ and 
 exchange the first and the third factors, then we obtain
 $s^{-1}_{\beta}(x'''_{(\beta^{-1})})\otimes\xi_{(\beta).i}x'_{(\beta)}\otimes\zeta_{(\beta^{-1}\alpha\beta).i}x''_{(\beta^{-1}\alpha\beta)}=s^{-1}_{\beta}(x'''_{(\beta^{-1})})\otimes x''_{(\beta).i}\xi_{(\beta).i}\otimes
  \varphi_{\beta^{-1}}(x'_{(\alpha)})\zeta_{(\beta^{-1}\alpha\beta).i}$.
 Applying $\mu_{\beta}\otimes T_{\beta^{-1}\alpha\beta}$ on both sides, 
 and observing that
 $s^{-1}_{\beta}(x'''_{(\beta^{-1})})x''_{(\beta)}\otimes\varphi_{\beta^{-1}}(x'_{(\alpha)})=
 \varepsilon(x'_{(\beta)})1_{\beta}\otimes\varphi_{\beta^{-1}}(x'_{(\alpha)})=1_{\beta}\otimes\varphi_{\beta^{-1}}(x)$,
 we obtain~\eqref{e:Asterix}.
\end{proof}

\begin{proof}[Proof of Theorem~\ref{thm:Asterix}]\label{p:proof-Asterix}
 By Theorem~\ref{thm:DH}, the \Tcoalg $D(H)$ defined as in
 Theorem~\ref{thm:dual} satisfies the three conditions 
 of Theorem~\ref{thm:Asterix}.
 We still have to check that these three conditions determinate 
 the \Tcoalg structure on $D(H)$, i.e., that the multiplication
 and the comultiplication we gave in the definition of $D(H)$ 
 are uniquely defined by the requirements of Theorem~\ref{thm:Asterix}. 
 
 \begin{sentence}{Product}
  Suppose that $D(H)$ satisfies the three conditions in 
  Theorem~\ref{thm:Asterix}.
 Let us check, that, given $\alpha\in\pi$, the multiplication on 
 $D_{\alpha}=H_{\alpha^{-1}}\circledast\Hdualcop$
 must satisfy the relations~\eqref{e:Golberg-Variationen}. 
 Let $p_{\alpha}$ be defined as in~\eqref{e:pro-Asterix}.
 The bijectivity of $p_{\alpha}$
 implies~\eqref{e:Golberg-Variationen-a}. 
 Let $h$ be in $\overline{H}_{\alpha}$ and let
 $f\in H^{\ast}_{\gamma}\subseteq\Hdualcop_{1}$, with $\gamma\in\pi$.
 By Lemma~\ref{e:Asterix}, we have
 \begin{equation}\label{e:Asterix-2}\begin{split}
  s^{-1}_{\gamma^{-1}}& \bigl((h\circledast\varepsilon)'''_{(\gamma)}\bigr)\,(e_{\gamma.i}\circledast\varepsilon)\,
  (h\circledast\varepsilon)'_{(\gamma^{-1})}\otimes(1_{\gamma\alpha^{-1}\gamma^{-1}}\circledast e^{\gamma.i})\,(h\circledast\varepsilon)''_{(\gamma\alpha\gamma^{-1})}\\
  & = (e_{\gamma.i}\circledast\varepsilon)\otimes\varphi_{\gamma}(h\circledast\varepsilon)\,(1_{\gamma\alpha\gamma^{-1}}\circledast e^{\gamma.i})\text{.}
 \end{split}\end{equation}
By observing that
$h'_{\overline{(\gamma^{-1})}}\otimes h''_{\overline{(\gamma\alpha\gamma^{-1})}}\otimes h'''_{\overline{(\gamma)}} =\varphi_{\gamma\alpha}(h'_{(\alpha^{-1}\gamma\alpha)})\otimes\varphi_{\gamma}(h'_{(\alpha^{-1})})\otimes h'''_{\gamma^{-1}}$,
if we compute the left-hand side in~\eqref{e:Asterix-2}, then we obtain
\begin{equation*}\begin{split}
 s^{-1}_{\gamma^{-1}}\bigl(& (h\circledast \varepsilon)'''_{(\gamma)}\bigr)\,(e_{\gamma.i}\circledast\varepsilon)\,
(h\circledast\varepsilon)'_{(\gamma^{-1})}\otimes(1_{\gamma\alpha^{-1}\gamma^{-1}}\circledast e^{\gamma.i})\,(h\circledast\varepsilon)''_{(\gamma\alpha\gamma^{-1})}\\ 
& =\bigl((\varphi_{\gamma^{-1}}\circ s_{\gamma})^{-1}(h'''_{(\gamma^{-1})})\circledast\varepsilon\bigr)
\,(e_{\gamma.i}\circledast\varepsilon)\,\bigl(\varphi_{\gamma\alpha}(h'_{(\alpha^{-1}\gamma\alpha)})\circledast\varepsilon\bigr)\otimes\\
&\qquad\qquad\otimes(1_{\gamma\alpha^{-1}\gamma^{-1}}\circledast e^{\gamma.i})\,\bigl(\varphi_{\gamma}(h'_{(\alpha^{-1})})\circledast\varepsilon\bigr)\\
& = \bigl((\varphi_{\gamma}\circ s^{-1}_{\gamma})(h'''_{(\gamma^{-1})})e_{\gamma.i}
\varphi_{\gamma\alpha}(h'_{(\alpha^{-1}\gamma\alpha)})\circledast\varepsilon\bigr)\otimes
\bigl(\varphi_{\gamma}(h'_{(\alpha^{-1})})\circledast e^{\gamma.i}\bigr)\text{,}
\end{split}\end{equation*}
where, in the last passage, we used both the fact that the 
immersion of $\overline{H}$ in $D(H)$ is a morphism of \Tcoalgs 
and~\eqref{e:Golberg-Variationen-a}.
So, if we apply $\varphi_{\gamma^{-1}}\otimes\varphi_{\gamma^{-1}}$ on both sides 
of~\eqref{e:Asterix-2}, then we get
\begin{multline}\label{e:Asterix-3}
 \bigl(s^{-1}_{\gamma}(h'''_{(\gamma^{-1})})e_{\gamma.i}\varphi_{\alpha}(h'_{(\alpha^{-1}\gamma\alpha)})\circledast\varepsilon\bigr)\otimes
 \bigl(h''_{(\alpha)}\circledast e^{\gamma.i}\bigr)\\
 =(e_{\gamma.i}\circledast\varepsilon)\otimes (h\circledast\varepsilon)(1_{\alpha^{-1}}\circledast e^{\gamma.i})\text{.}
\end{multline}
If we evaluate both terms of~\eqref{e:Asterix-3} against 
$\langle f\otimes 1_{1},\_\rangle\otimes D_{\alpha}(H)$, then on the left-hand side we get
$\bigl\langle f, s^{-1}_{\gamma}e_{\gamma.i}\varphi_{\alpha}(h'_{(\alpha^{-1}\gamma\alpha)})\bigr\rangle h''_{(\alpha)}\circledast e^{\gamma.i}
 = h''_{(\alpha)}\circledast \bigl\langle f, s^{-1}_{\gamma}\_\varphi_{\alpha}(h'_{(\alpha^{-1}\gamma\alpha)})\bigr\rangle$,
while on the right-hand side we get
$\langle f, e_{\gamma.i}\rangle (h\circledast\varepsilon)\,(1_{\alpha^{-1}}\circledast 
 e^{\gamma.i})=(h\circledast\varepsilon)\,(1_{\alpha^{-1}}\circledast f)$.
\end{sentence}

\begin{sentence}{Comultiplication}
Let us check that the comultiplication on $D(H)$ is also unique. 
Given $\alpha,\beta,\gamma\in\pi$, $h\in \overline{H}_{\alpha\beta}$, and $f\in H^{\ast}_{\gamma}$, we have
%$\Delta_{\alpha,\beta}(h\circledast f) = \Delta_{\alpha,\beta}\bigl((1_{\alpha^{-1}}\circledast f)\,(h\circledast\varepsilon)\bigr)
% = (1_{\alpha^{-1}}\circledast f')\, (h'_{\overline{(\alpha)}}\circledast\varepsilon)\otimes
% (1_{\alpha^{-1}}\circledast f'')\,(h''_{\overline{(\beta)}}\circledast\varepsilon)
% = \varphi_{\beta}\bigl(h'_{(\beta^{-1}\alpha^{-1}\beta)}\bigr)\circledast f'\otimes h''_{(\beta^{-1})}\circledast f''$.
\begin{equation*}\begin{split}
\Delta_{\alpha,\beta}(h\circledast f) & = \Delta_{\alpha,\beta}\bigl((1_{\alpha^{-1}}\circledast f)\,(h\circledast\varepsilon)\bigr)
 \\ & = (1_{\alpha^{-1}}\circledast f')\, (h'_{\overline{(\alpha)}}\circledast\varepsilon)\otimes
 (1_{\alpha^{-1}}\circledast f'')\,(h''_{\overline{(\beta)}}\circledast\varepsilon)
 \\ & = \Bigl(\varphi_{\beta}\bigl(h'_{(\beta^{-1}\alpha^{-1}\beta)}\bigr)\circledast 
        f'\Bigr)\otimes\Bigl( h''_{(\beta^{-1})}\circledast f''\Bigr)\text{.}
\end{split}\end{equation*}
\end{sentence}
\end{proof}

\begin{rmk}
 The quantum double of a Hopf algebra can be obtained also 
 via the Majid bicrossproduct~\cite{Majid-Phys}.
 This is true also in the crossed case~\cite{Zunino-bicross}. 
 Recently, Vainerman~\cite{Vaynermann-HQFT} constructed some
 examples of nontrivial \Tcoalgs with non-isomorphic components.
 These examples can also be interpreted 
 ---and, eventually, enlarged--- as bicrossproducts of \Tcoalgs.
\end{rmk}

\begin{rmk}
 Street~\cite{Street-double} has proved that, starting from any 
 (not necessarily finite-dimensional) Hopf algebra $H$, 
 it is possible to construct, via Tannaka Theory, a coquasitriangular 
 Hopf algebra $D^{\ast}(H)$ such that, when $H$ is finite-dimensional, 
 $D^{\ast}(H)=\bigl(D(H)\bigr)^{\ast}$. Tannaka Theory for \Talgs was
 developed by the author in~\cite{Zunino-Tannaka} where, contextually, 
 it is provided an analog for the co-double construction~\cite{ReSe} 
 in the case of a \Talg.
\end{rmk}

\section{Ribbon \protect\Tcoalgs}\label{fur-Alina-2}
 The notion of a ribbon Hopf algebra~\cite{RT} is 
 generalized to the case of a \Tcoalg in~\cite{Tur-CPC}.
 
 \paragraph{\scshape Drinfeld elements}
  Following~\cite{Virelizi}, we set
 $u_{\alpha}\eqdef (s_{\alpha^{-1}}\circ\varphi_{\alpha})(\zeta_{(\alpha^{-1}).i})\xi_{(\alpha).i}$
 and $u\eqdef\{u_{\alpha}\}_{\alpha\in\pi}$. 
 The $u_{\alpha}$ are called \textit{Drinfeld elements} of $H$.
 When $\pi=\{1\}$, we recover the usual definition 
 of Drinfeld element of a quasitriangular Hopf algebra.
 
 The following properties of $u$ are proved in~\cite{Virelizi}. 
 Let $\alpha$ and $\beta$ be in $\pi$ and let $h$ be in $H_{\alpha}$.
 \begin{eqnenumerate}\label{Ovid}
  \item $u_{1}=s_{1}(\zeta_{(1).i})\xi_{(1).i}$.\label{Ovid-1}
  \item $u_{\alpha}$ is invertible with inverse 
        $u^{-1}_{\alpha}=s^{-1}_{\alpha}\Bigl(\tilde{\zeta}_{(\alpha^{-1}).i}\Bigr)\tilde{\xi}_{(\alpha).i}$.
        Moreover we have\label{Ovid-2} 
        $u_{\alpha}^{-1}  =(s_{\alpha}^{-1}\circ s^{-1}_{\alpha^{-1}})(\zeta_{(\alpha).i})\xi_{(\alpha).i}=
         =\xi_{(\alpha).i}(s_{\alpha^{-1}}\circ s_{\alpha})(\zeta_{(\alpha).i})$.
  \item $(u_{\alpha\beta})'_{(\alpha)}\otimes (u_{\alpha\beta})''_{(\beta)}=
         \tilde{\xi}_{(\alpha).i}\tilde{\zeta}_{(\alpha).j}u_{\alpha}\otimes\tilde{\zeta}_{(\beta).i}
         \varphi_{\alpha^{-1}}(\tilde{\xi}_{(\alpha\beta\alpha^{-1}).j})u_{\beta}$.\label{Ovid-3}
  \item $\varepsilon(u_{1})=1$.\label{Ovid-4}
  \item $s_{\alpha^{-1}}(u_{\alpha^{-1}})u_{\alpha}=u_{\alpha}s_{\alpha^{-1}}(u_{\alpha^{-1}})$.\label{Ovid-5}
  \item $\varphi_{\beta}(u_{\alpha})=u_{\beta\alpha\beta^{-1}}$.\label{Ovid-6}
  \item $(s_{\alpha^{-1}}\circ s_{\alpha}\circ\varphi_{\alpha})(h)=u_{\alpha}hu_{\alpha}^{-1}$.\label{Ovid-7}
  \item $u_{\alpha}s_{\alpha^{-1}}(u_{\alpha^{-1}})h=
         \varphi_{\alpha^{2}}(h)u_{\alpha}s_{\alpha^{-1}}(u_{\alpha^{-1}})$.\label{Ovid-8}
 \end{eqnenumerate}

\noindent Notice that, by~\eqref{Ovid-7}, we have
$(s_{\alpha^{-1}}\circ s_{\alpha})(h) =\varphi_{\alpha^{-1}}(u_{\alpha})\varphi_{\alpha^{-1}}(h)\varphi_{\alpha^{-1}}(u^{-1}_{\alpha})=u_{\alpha}\varphi_{\alpha^{-1}}(h)u^{-1}_{\alpha}$.
In particular, for $h=u_{\alpha}$ we obtain
\begin{eqnenumerate}\addtocounter{equation}{-1}\setcounter{enumi}{8}
  \item $(s_{\alpha^{-1}}\circ s_{\alpha})(u_{\alpha})=u_{\alpha}$.\label{Ovid-9}
\end{eqnenumerate}

\paragraph{\scshape Definition of a ribbon \protect\Tcoalg (first version)}
 Let $H$\label{p:DrinfeldElements} 
 be a quasitriangular \Tcoalg.
According to~\cite{Tur-CPC}, we say that $H$ is a 
\textit{ribbon \Tcoalg} if it is endowed with a family
 $\theta=\{\theta_{\alpha}\vert\theta_{\alpha}\in H_{\alpha}\}_{\alpha\in\pi}$, the \textit{twist,}
 such that $\theta_{\alpha}$ is invertible for any $\alpha\in\pi$ and 
 the following conditions are satisfied
 for any $\alpha,\beta\in\pi$ and $h\in H_{\alpha}$.
 \begin{eqnenumerate}\label{bkdingnum}
  \item $\varphi_{\alpha}(h) = \theta^{-1}_{\alpha}h\theta_{\alpha}$.\label{bkdingnumeral1}
  \item $s_{\alpha}(\theta_{\alpha})=\theta_{\alpha^{-1}}$.\label{bkdingnumeral2}
  \item $(\theta_{\alpha\beta})'_{(\alpha)}\otimes (\theta_{\alpha\beta})''_{(\beta)}=
        \theta_{\alpha}\zeta_{(\alpha).i}\xi_{(\alpha).j}\otimes \theta_{\beta}        
        \varphi_{\alpha^{-1}}(\xi_{(\alpha\beta\alpha^{-1}).i})\zeta_{(\beta).j}$.\label{p:bkdin3}\label{bkdingnumeral3}
  \item $\varphi_{\beta}(\theta_{\alpha})=\theta_{\beta\alpha\beta^{-1}}$.\label{bkdingnumeral4}
\end{eqnenumerate}

\noindent Notice that $(H_{1},R_{1,1},\theta_{1})$ is a ribbon Hopf algebra 
in the usual sense.

If $H=(H,R,\theta)$ is a ribbon \Tcoalg, then, for any $\alpha\in\pi$, 
we obtain the following properties.
 \begin{eqnenumerate}\label{Silence}
  \item $\varphi_{\alpha^{-1}}(h)=\theta_{\alpha}h\theta_{\alpha}^{-1}$ for any $h\in H_{\alpha}$ 
        (this follows by~\ref{bkdingnumeral1}).\label{Silence-1}
  \item $\varepsilon(\theta_{1}) = 1$ 
        (this is because $H_{1}$ is a ribbon Hopf algebra).\label{Silence-2}
  \item $\theta_{1}$ is central (by the same reason).\label{Silence-3}
  \item $\theta_{\alpha}\varphi_{\alpha}(h)=h\theta_{\alpha}$ for any $h\in H_{\alpha}$ 
        (this follows by~\ref{bkdingnumeral1}).\label{Silence-4}
 \end{eqnenumerate}    

 \noindent Moreover, it is proved in~\cite{Virelizi} that we have

 \begin{eqnenumerate}\addtocounter{equation}{-1}\setcounter{enumi}{4}
  \item $\theta_{\alpha}u_{\alpha}=u_{\alpha}\theta_{\alpha}$\label{Silence-5} and
  \item $\theta^{-2}_{\alpha}=s_{\alpha^{-1}}(u_{\alpha^{-1}})u_{\alpha}=
         u_{\alpha}s_{\alpha^{-1}}(u_{\alpha^{-1}})$.\label{Silence-6}\label{p:Amen}
 \end{eqnenumerate}

If for any $\alpha,\beta\in\pi$, we introduce the notation 
 \begin{subequations}\label{e:continuum}
 \begin{multline}\label{e:continuum-a}
  Q_{\alpha,\beta} =\bigl(\sigma\circ(\varphi_{\alpha^{-1}}\otimes H_{\alpha})\bigr)(R_{\alpha\beta\alpha^{-1},\alpha})R_{\alpha,\beta}=
  \\
   =\zeta_{(\alpha).i}\xi_{(\alpha).j}\otimes\varphi_{\alpha^{-1}}(\xi_{(\alpha\beta\alpha^{-1}).i})\zeta_{(\beta).j}\text{,}
 \end{multline}
 then~\eqref{bkdingnumeral3} can be written in the form
 $(\theta_{\alpha\beta})'_{(\alpha)}\otimes (\theta_{\alpha\beta})''_{(\beta)}=(\theta_{\alpha}\otimes\theta_{\beta})Q_{\alpha,\beta}$.
 Moreover, if we set
 $\tilde{Q}_{\alpha,\beta}\eqdef Q_{\alpha,\beta}^{-1}=
  \tilde{\xi}_{(\alpha).i}\tilde{\zeta}_{(\alpha).j}\otimes
  \tilde{\zeta}_{(\beta).i}\varphi_{\alpha^{-1}}(\tilde{\xi}_{(\alpha\beta\alpha^{-1}).j})$,
 we can rewrite~\eqref{Ovid-3} in the form
 \begin{equation}\label{e:delta-u}
  (u_{\alpha\beta})'_{(\alpha)}\otimes (u_{\alpha\beta})''_{(\beta)}=\tilde{Q}_{\alpha,\beta}(u_{\alpha}\otimes u_{\beta})\text{.}
 \end{equation}
 We observe that the conjugation preserves $Q$, i.e., 
 that, for any $\alpha,\beta,\gamma\in\pi$, we have
 $(\varphi_{\alpha}\otimes\varphi_{\alpha})(Q_{\beta,\gamma}) = Q_{\alpha\beta\alpha^{-1},\alpha\gamma\alpha^{-1}}$.
 \end{subequations}

 \paragraph{\scshape Definition of a ribbon \protect\Tcoalg 
 (second version)}
   We\label{p:second-def} can define a ribbon \Tcoalg 
   in an equivalent way as a quasitriangular \Tcoalg $H$ endowed 
   with a family
 $v=\{v_{\alpha}\vert v_{\alpha}\in H_{\alpha}\}_{\alpha\in\pi}$ satisfying the following conditions.
 \begin{eqnenumerate}\label{e:martinello}
  \item $hv_{\alpha}=v_{\alpha}\varphi_{\alpha^{-1}}(h)$ for any $h\in H_{\alpha}$.\label{dingnumeral1}
  \item $v^{2}_{\alpha}=u_{\alpha}s_{\alpha^{-1}}(u_{\alpha^{-1}})$.\label{dingnumeral2}
  \item $(v_{\alpha\beta})'_{(\alpha)}\otimes(v_{\alpha\beta})''_{(\beta)}= 
        \tilde{Q}_{\alpha,\beta}(v_{\alpha}\otimes v_{\beta})$.\label{dingnumeral3}
  \item $s_{\alpha}(v_{\alpha})=v_{\alpha^{-1}}$.\label{dingnumeral4}
  \item $\varphi_{\beta}(v_{\alpha})=v_{\beta\alpha\beta^{-1}}$.\label{dingnumeral5}
 \end{eqnenumerate}

 The proof of the equivalence of the two definitions can be easily 
 obtained by modifying the analogous easy proof for Hopf algebras.

 \section{The quantum double of a semisimple 
          \protect\Tcoalg}\label{fur-Alina-3}
  The quantum double of a 
  semisimple Hopf algebra over a field of characteristic zero 
  is both semisimple and modular (see~\cite{EG}).
   We start this section by recalling the definition of a 
  semisimple \Tcoalg~\cite{Virelizi}, a modular Hopf algebra~\cite{RT1}, 
  and a modular \Tcoalg~\cite{Tur-CPC,Virelizi}.
  After that, given any totally-finite \Tcoalg $H$, we discuss 
  the relation between $D(H)$ and the quantum double of $H_{\pk}$.
  Finally, we discuss the semisimplicity and the modularity 
  of the quantum double $D(H)$ of a semisimple \Tcoalg 
  $H$ over a field of characteristic zero. 
  In particular, we prove that $D(H)$ is semisimple 
  if and only if $H$ is totally-finite.
  Moreover, when $H$ is totally-finite, $D(H)$ is also modular.

 \paragraph{\scshape Basic definitions}
  Let $H$ be a \Tcoalg. 
  We say that $H$ is \textit{semisimple} when any algebra $H_{\alpha}$ 
  (with $\alpha\in\pi$) is semisimple. 
  It is proved in~\cite{Virelizi} that $H$ is semisimple 
  if and only if $H_{1}$ is semisimple.
  Further, by~\cite{Sweedler-integral}, 
  infinite-dimensional Hopf algebras over a field are never semisimple. 
  It follows that a necessary condition for $H$ to be semisimple 
  is that $H_{1}$ is finite-dimensional.

  Let $H_{1}=(H_{1},R_{1}=\xi_{1.i}\otimes\zeta_{1.i},\theta_{1})$ 
  be a ribbon Hopf algebra. 
 Given a finite-dimensional representation $V$ of $H_{1}$, 
 and a \Hlinear endomorphism $f\colon V\to V$, the 
 \textit{quantum trace $\trq(f)$ of $f$} is defined as
 $\trq(f)\eqdef\tr(u_{1}\theta_{1} f)$,
 where $u_{1}=s_{1}(\zeta_{1.i})\xi_{1.i}$ and $\tr(\cdot)$ is 
 the usual trace of endomorphisms.
 $V$ is said to be \textit{negligible} when $\trq(\Id_{V})=0$.

 A \textit{modular Hopf algebra} $H_{1}$ is a ribbon Hopf algebra 
 endowed with a finite family of simple finite-dimensional \Honemodules 
 $\{V_{i}\}_{i\in I}$ satisfying the following conditions.
 \begin{itemize}
  \item There exists $0\in I$ such that $V_{0}=\Bbbk$ 
        (with the structure of \Honemodule 
        given by the comultiplication).
  \item For any $i\in I$, there exists $i^{\ast}\in I$ such that $V_{i^{\ast}}$ 
        is isomorphic to $V_{i}^{\ast}$.
  \item For any $j,k\in I$, the \Honemodule $V_{j}\otimes V_{k}$ 
        is isomorphic to a finite sum of certain elements of $\{V_{i}\}_{i\in I}$,
        possibly with repetitions, and a negligible \Honemodule.
  \item Let $(S_{i,j})_{i,j\in I}$ denote the square matrix 
        whose entry $S_{i,j}$ is the quantum trace of
        the endomorphism $x\mapsto \zeta_{1.k}\xi_{1.l}\otimes\xi_{1.k}\zeta_{1.l}x$ of
        $V_{i}\otimes V_{j^{\ast}}$. Then $S[H]$ is invertible.
 \end{itemize}

A \textit{modular \Tcoalg}~\cite{Tur-CPC} 
is a ribbon \Tcoalg $H$ such that its component $H_{1}$ 
is a modular Hopf algebra.

 \begin{thm}\label{thm:ssmod} 
  The quantum double $D(H)$ of a \Tcoalg $H$ over a
  field of characteristic $0$ is semisimple if and only if
  $H$ is totally-finite, and, in that case, $D(H)$ is also modular.
 \end{thm}
 
To prove Theorem~\ref{thm:ssmod} we firstly need to discuss 
how the quantum double of a totally-finite
\Tcoalg $H$ and the quantum double of the packed form of $H$ are
related.

\paragraph{\scshape 
 The structure of a totally finite \protect\Tcoalg}\label{par:tf-tf}
 When $H$ is a totally-finite \Tcoalg, the Hopf algebra 
 $H^{\innersym}_{1}=(H^{\outersym})_{\pk}$ is the dual of a certain Hopf
 algebra $H_{\pk}$.
 An easy computation shows that $H_{\pk}$ satisfies the 
 following conditions.
 \begin{itemize}
  \item As an algebra, $H_{\pk}=\bigoplus_{\alpha\in\pi}H_{\alpha}$.
  \item The comultiplication $\Delta_{\pk}$ is obtained setting
        $\Delta_{\pk}\eqdef\sum_{\beta,\gamma\colon \beta\gamma=\alpha}\Delta_{\beta,\gamma}(h)$,
        for any $h\in H_{\alpha}\subset H_{\pk}$.
  \item The counit $\varepsilon_{\pk}$ is given by
        $\varepsilon_{\pk}{\raisebox{-.5ex}{$\vert_{H_{1}}$}}= \varepsilon$
        and $\varepsilon_{\pk}{\raisebox{-.5ex}{$\vert_{H_{\alpha}}$}}= 0$, for
        any $\alpha\in\pi$, $\alpha\neq1$.
  \item The antipode is given by $s_{\pk}=\sum_{\alpha\in\pi}s_{\alpha}$.
  \end{itemize}
 
\noindent Moreover, we have a group homomorphism 
    $\map{\varphi_{\pk}}{\pi}{\!\!\!\!\!\Aut(H_{\pk})}{\alpha}{\bigoplus_{\beta\in\pi}\varphi_{\beta}^{\alpha}}$.
 
Conversely, let $H_{\grande}$ be a finite-dimensional Hopf algebra 
(with antipode $s_{\grande}$, counit $\varepsilon$ and comultiplication 
$\Delta_{\grande}\/$), endowed with a family of subcoalgebras $\{H_{\alpha}\}_{\alpha\in\pi}$ 
and a group homomorphism 
$\map{\varphi_{\grande}}{\pi}{\Aut(H_{\grande})}{\alpha}{\displaystyle\varphi_{\grande,\alpha}}$,
such that the following conditions hold.
 \begin{itemize}
  \item $H_{\grande}$ is, as an algebra, the product $H_{\alpha}$ ($\alpha\in\pi\/$).
  \item For any $\alpha\in\pi$, we have 
        $\Delta_{\grande}(H_{\alpha})\subset\bigoplus_{\beta,\gamma\text{ s.t. }\beta\gamma=\alpha}(H_{\beta}\otimes H_{\gamma})$.
  \item For any $\alpha\in\pi\setminus\{1\}$, $H_{\alpha}\subset\Ker\,\varepsilon$.
  \item For any $\alpha\in\pi$, $s_{\grande}(H_{\alpha})= H_{\alpha^{-1}}$.
  \item For any $\alpha,\beta\in\pi$, the image of $H_{\alpha}$ under $\varphi_{\grande,\beta}$ 
        lies in $H_{\beta\alpha\beta^{-1}}$.
 \end{itemize}
 
\noindent To $H_{\grande}$ corresponds, 
in the obvious way, a \Tcoalg $H$ such that $H_{\pk}=H_{\grande}$.
In particular, for any $\alpha,\beta\in\pi$, the component $\Delta_{\alpha,\beta}\colon H_{\alpha\beta}\to H_{\alpha}\otimes H_{\beta}$ 
of the comultiplication is given by
$
 H_{\alpha\beta}\hookrightarrow H_{\grande}\xrightarrow{\Delta_{\grande}}H_{\grande}\otimes H_{\grande}
\xrightarrow{p_{\alpha}\otimes p_{\beta}}H_{\alpha}\otimes H_{\beta}$,
where $p_{\alpha}$ and $p_{\beta}$ are the canonical projections 
of $H_{\grande}$ on $H_{\alpha}$ and, respectively, $H_{\beta}$.
 
\paragraph{\scshape The quantum double 
           of a totally-finite \protect\Tcoalg}
 Let $H$ be a totally-finite \Tcoalg. 
 We have seen that $D(H)$ is totally finite.
 We can construct in the usual way 
 the quantum double $D(H_{\pk})$ of $H_{\pk}$. 
 Neither as an algebra nor as a coalgebra 
 $D(H_{\pk})$ is isomorphic to $\bigl(D(H)\bigr)_{\pk}$
 since neither the multiplication nor the comultiplication 
 of $D(H_{\pk})$ depend on the conjugation.
  
 Now, the canonical embedding of vector spaces 
$H_{1}\otimes \Hdualcop_{1}\hookrightarrow \bigoplus_{\alpha\in\pi}H_{\alpha}\otimes\Hdualcop_{1}$ 
provides an embedding of Hopf algebras of $D_{1}(H)\hookrightarrow D(H_{\pk})$, 
so we can identify $D_{1}(H)$ with its image in $D(H_{\pk})$.
 Moreover, even if the universal \Rmatrix
 $R_{1,1}=(e_{1.i}\otimes\varepsilon)\otimes (1\otimes e^{1.i})$ of $D_{1}(H)$ and the universal 
 \Rmatrix $R_{\pk}=\sum_{\alpha\in\pi}R_{\alpha,\alpha}$ of $D(H_{\pk})$
 are different, for any $x\in H_{1}\otimes H^{\ast}$ we have
 \begin{equation}\label{e:Santiango-3}
   xR_{\pk}=xR_{1,1}\qquad\text{and}\qquad R_{\pk}x=R_{1,1}x\text{.}
 \end{equation}
 
 \paragraph{\scshape Factorizable \protect\Tcoalgs}
 We recall that a finite dimensional quasitriangular Hopf algebra
 $H$ with universal \Rmatrix $R=\xi_i\otimes\zeta_{i}$ is \textit{factorizable}
 (see~\cite{ReSe}) if the map 
 $\map{\lambda}{H^{\ast}}{H}{f}{\langle f,\zeta_{i}\xi_{j}\rangle \xi_{i}\zeta_{j}}$ is
 bijective. In particular, the quantum double of any finite-dimensional
 Hopf algebra is factorizable. It is proved in~\cite{Radford} that any
 factorizable Hopf algebra is unimodular.
  
 \begin{lemma}\label{l:factorizable}
  Let $H$ be a finite-type \Tcoalg. The Hopf algebra
  $\bigl(D_{1}(H),\xi_{(1).i}\otimes \zeta_{(1).i}\bigr)$ is factorizable.
 \end{lemma}
 
 \begin{proof}
  We can identify $D_1(H)$ as a subspace of $D(H_{\pk})$ and 
  $\bigl(D_1(H)\bigr)^{\ast}$ as a subspace of $\bigl(D(H_{\pk})\bigr)^{\ast}$.
  We only need to check that 
  $\lambda_{D_1(H)}=\lambda_{D(H_{\pk})}\vert_{\bigl(D_1(H)\bigr)^{\ast}}$.
  Let $R=\xi_i\otimes\zeta_i$ be the universal \Rmatrix of $D(H_{\pk})$. For any
  $f\in\bigl(D_1(H)\bigr)^{\ast}$ we have
  \begin{equation*}
   \lambda_{D(H_{\pk})}(f)=\langle f,\zeta_{i}\xi_{j}\rangle \xi_{i}\zeta_{j}=
   \langle f,\zeta_{(1).i}\xi_{(1).j}\rangle \xi_{(1).i}\zeta_{(1).j}=\lambda_{D_1(H)}(f)\text{.}
  \end{equation*}
 \end{proof}

 \paragraph{\scshape The quantum double of a semisimple \protect\Tcoalg}
 Let us consider the case of a semisimple 
 \Tcoalg $H$ over a field $\Bbbk$ of characteristic zero. 
 It was proved in~\cite{Virelizi} that,
 for any $\alpha\in\pi$,
 \begin{equation}\label{e:Mahler}
  s_{\alpha^{-1}}\circ s_{\alpha}=\Id_{H_{\alpha}}\text{.}
 \end{equation}
 
  \begin{lemma}\label{l:ref-a-3} 
  If $H$ is quasitriangular, then it is also ribbon by setting
  $\theta_{\alpha}=u_{\alpha}^{-1}$, for any $\alpha\in\pi$.
 \end{lemma}
 
 \begin{proof}
  Axiom~\ref{dingnumeral1} follows by~\eqref{Ovid-7} 
  and~\eqref{e:Mahler}. 
  Axiom~\ref{dingnumeral3} follows by~\eqref{e:delta-u}.
  Axiom~\ref{dingnumeral4} follows by~\eqref{Ovid-6}.
  Axiom~\ref{dingnumeral5} can be rewritten
  \begin{equation}\label{e:ref-a-3}
   s_{\alpha^{-1}}(u_{\alpha^{-1}}) = u_{\alpha^{-1}}\text{.}
  \end{equation}
  This follow by~\cite[Theorem~6(b)]{Virelizi}, by observing that, 
  in that formula,
  $g_{\alpha}=1_{\alpha}$ (by~\cite[Corollary~7]{Virelizi})
  $\hat{\varphi}(\alpha)=1$ (by~\cite[Theorem~7]{Virelizi}), 
  and $h_{\alpha}=1_{\alpha}$ (by~\cite[Lemma~16]{Virelizi},
  since the distinguished group-like element of $H^{\ast}_{1}$
  is equal to $\varepsilon$ because $H^{\ast}_{1}$ is semisimple by~\cite{LR1}).
  Finally axiom~\ref{dingnumeral2} follows by~\eqref{e:ref-a-3}.
 \end{proof}

 \begin{proof}[Proof of Theorem~\ref{thm:ssmod}]
  If $H$ is not totally-finite, $D_{1}(H)$ is not finite-dimensional 
  and so it is not semisimple.
 
  Suppose that $H$ is totally-finite. 
  $H_{\alpha}$ is semisimple for any $\alpha\in\pi$, therefore $H_{\pk}$ is semisimple. 
  It follows that $D(H_{\pk})$ is semisimple (see~\cite{Radford}). 
  Since $D_{1}(H)$ can be 
  identified with a subalgebra of $D(H_{\pk})$, also $D_{1}(H)$ is 
  semisimple, hence $D(H)$ is semisimple. By Lemma~\ref{l:ref-a-3}, 
  $D(H)$ has a natural structure of a ribbon \Tcoalg.
  
  The rest of the proof follows the streamline of the proof of Lemma~1.1
  in~\cite{EG} to which we refer for further details.
  Let $C\bigl(D_1(H)\bigr)\subset\bigl(D_1(H)\bigr)^{\ast}$ be the ring of 
  characters of $D_1(H)$ and let $Z\bigl(D_1(H)\bigr)$ be the
  center of $D_1(H)$. By Lemma~\ref{l:factorizable}, the restriction
  $\Lambda=\lambda\vert_{C\bigl(D_1(H)\bigr)}\colon C\bigl(D_1(H)\bigr)\to Z\bigl(D_1(H)\bigr)$
  is an isomorphism of algebras.
  Let $\Irr\bigl(D_{1}(H)\bigr)=\{V_{i}\vert 0\leq i\leq m\}$ be a set 
  of representatives for the isomorphism classes of the irreducible 
  representations of $D(H_{1})$ such that $V_{0}=\Bbbk$. Then
  $B=\{\chi_{j^{\ast}}=\tr\vert_{V_{j^{\ast}}}\colon 0\leq j\leq m\}$ is a linear basis of
  $C\bigl(D(H_1)\bigr)$. We also observe that the set $C=\{e_j\colon 0\leq j\leq m\}$ 
  of central primitive idempotents of $D_{1}(H)$ is a basis of
  $Z\bigl(D(H_1)\bigr)$. Now, $S[D_1(H)]=DA$, where $A$ is the
  invertible matrix which represents $\Lambda$ with respect to
  the bases $B$ and $C$, while $D=\mbox{diag}(\dim V_i)$.
\end{proof}

\section{The ribbon extension of a
 quasitriangular \protect\Tcoalg}\label{sec:ribbonificator}
 Let $H$ be any quasitriangular \Tcoalg 
 (not necessarily of finite-type).
 We describe how to obtain a ribbon \Tcoalg $\RT(H)$ such that, 
 when $\pi=\{1\}$, we recover the construction described in~\cite{RT}.

\paragraph{\scshape Definition of $\RT(H)$} 
 The \textit{ribbon extension} 
 $RT(H)$ of a quasitriangular \Tcoalg $H$ 
 is the \Tcoalg defined as follows. 
 \begin{itemize}
  \item For any $\alpha\in\pi$, the \alphath component of $RT(H)$, 
        denoted $RT_{\alpha}(H)$, is the vector space whose elements 
        are formal expressions $h+kv_{\alpha}$, with $h, k\in H_{\alpha}$. 
        The sum is given by
        $(h+kv_{\alpha}) + (h'+k'v_{\alpha})\eqdef (h + h') + (k + k')v_{\alpha}$,
        for any $h, h', k, k'\in H_{\alpha}$. 
        The multiplication is obtained by requiring 
        $v^{2}_{\alpha}=u_{\alpha}s_{\alpha^{-1}}(u_{\alpha^{-1}})$, i.e., by setting, 
        for any $h, h', k, k'\in H_{\alpha}$,
        \begin{equation*}
         (h+kv_{\alpha})\,(h'+k'v_{\alpha})= 
         \bigl(hh'+k\varphi_{\alpha}(k')u_{\alpha}s_{\alpha^{-1}}(u_{\alpha^{-1}})\bigr)+
           \bigl(hk'+k\varphi_{\alpha}(k')\bigr)v_{\alpha}\text{.}
        \end{equation*}

        We identify $H_{\alpha}$ with the subset $\{h+0v_{\alpha}\vert h\in H_{\alpha}\}$ 
        of $RT_{\alpha}(H)$. 
        The algebra $RT_{\alpha}(H)$ is unitary with unit $1_{\alpha}= 1_{\alpha}+0v_{\alpha}$. 
        Moreover, for any $\alpha,\beta\in\pi$, 
        we have $R_{\alpha,\beta}\in H_{\alpha}\otimes H_{\beta}\subset RT_{\alpha}(H)\otimes RT_{\beta}(H)$.
  \item The comultiplication is given by
        \begin{equation*}
        \Delta_{\alpha,\beta} (h+kv_{\alpha\beta})
         =\Delta_{\alpha,\beta}(h)+\Delta_{\alpha,\beta}(k)\tilde{Q}_{\alpha,\beta}(v_{\alpha}\otimes v_{\beta})\text{,}
        \end{equation*}
        for any $h, k\in H_{\alpha}$ and $\alpha,\beta\in\pi$.
        The counit is given by
        $\langle\varepsilon, h+kv_{\alpha}\rangle =\langle\varepsilon, h\rangle+\langle\varepsilon,k\rangle$,
        for any $h, k\in H_{1}$.
  \item The antipode is given by
        $s_{\alpha}(h+kv_{\alpha}) = s_{\alpha}(h)+(s_{\alpha}\circ\varphi_{\alpha^{-1}})(k)v_{\alpha^{-1}}$,
        for any $h, k\in H_{\alpha}$ and $\alpha\in\pi$.
  \item Finally, the conjugation is given by
        $\varphi_{\beta}(h+kv_{\alpha})=\varphi_{\beta}(h)+\varphi_{\beta}(k)v_{\beta\alpha\beta^{-1}}$,
        for any $h, k\in H_{\alpha}$ and $\alpha,\beta\in\pi$.
 \end{itemize}
 
 \begin{thm}\label{thm:RT}
  $\RT(H)$ is a ribbon \Tcoalg.
 \end{thm}

\begin{proof}
 Since computations involved to prove in detail Theorem~\ref{thm:RT}
 are relatively long, we only prove the coassociativity of $\Delta$,
 while the rest is left to the reader. To prove the
 multiplicativity of $\Delta$ you need to check before the
 relation
 \begin{equation*}
  \Delta_{\alpha,\beta} \bigl(u_{\alpha\beta}s_{(\alpha\beta)^{-1}}(u_{(\alpha\beta)^{-1}})\bigr)=
  \tilde{Q}_{\alpha,\beta}(\varphi_{\alpha}\otimes\varphi_{\beta})(\tilde{Q}_{\alpha,\beta})
  u_{\alpha}s_{\alpha^{-1}}(u_{\alpha^{-1}})\otimes u_{\beta}s_{\beta^{-1}}(u_{\beta^{-1}})\text{.}
 \end{equation*}
 (for any $\alpha,\beta\in\pi$\/).
 To prove that $s_{\RT(H)}$ is an antipode, 
 one first needs to show that it is antimultiplicative.
 
\begin{sentence}{Coassociativity}
 We need to check that, for any $h,k\in H_{\alpha\beta\gamma}$, with $\alpha,\beta,\gamma\in\pi$, we have
 \begin{multline}\label{e:Rossini}
  \bigl((\Delta_{\alpha,\beta}\otimes \RT_{\gamma}(H))\circ\Delta_{\alpha\beta,\gamma}\bigr)(h+kv_{\alpha\beta\gamma})\\ =
  \bigl((\RT_{\alpha}(H)\otimes\Delta_{\beta,\gamma})\circ \Delta_{\alpha,\beta\gamma}\bigr)(h+kv_{\alpha\beta\gamma})\text{.}
 \end{multline}
 By the linearity and the multiplicativity of the comultiplication 
 in $\RT(H)$ and the coassociativity of the comultiplication in $H$,
 by computing the left-hand side of~\eqref{e:Rossini}, we obtain
 $\bigl(( \Delta_{\alpha,\beta}\otimes \RT_{\gamma}(H))\circ\Delta_{\alpha\beta,\gamma}\bigr)(h+kv_{\alpha\beta\gamma})
   = h'_{(\alpha)}\otimes h''_{(\beta)}\otimes h'''_{(\gamma)} +
  (k'_{(\alpha)}\otimes k''_{(\beta)}\otimes k'''_{(\gamma)})
  \bigl((\Delta_{\alpha,\beta}\otimes \RT_{\gamma}(H))\circ\Delta_{\alpha\beta,\gamma}\bigr)(v_{\alpha\beta\gamma})$,
 while, by computing the right-hand, side we obtain
 $\bigl((\RT_{\alpha}(H)\otimes \Delta_{\beta,\gamma})\circ \Delta_{\alpha,\beta\gamma}\bigr)(h+kv_{\alpha\beta\gamma})
  = h'_{(\alpha)}\otimes h''_{(\beta)}\otimes h'''_{(\gamma)} +
  (k'_{(\alpha)}\otimes k''_{(\beta)}\otimes k'''_{(\gamma)})
  \bigl((\RT_{\alpha}(H)\otimes \Delta_{\beta,\gamma})\circ \Delta_{\alpha,\beta\gamma}\bigr)(v_{\alpha\beta\gamma})$,
 so we only need to check that we have
 \begin{equation}\label{e:Roma}
  \bigl((\Delta_{\alpha,\beta}\otimes \RT_{\gamma}(H))\circ\Delta_{\alpha\beta,\gamma}\bigr)(v_{\alpha\beta\gamma})=
  \bigl((\RT_{\alpha}(H)\otimes \Delta_{\beta,\gamma})\circ \Delta_{\alpha,\beta\gamma}\bigr)(v_{\alpha\beta\gamma})\text{.}
 \end{equation}
 By computing the left-hand side of~\eqref{e:Roma}, we get
 $
 \bigl((\Delta_{\alpha,\beta} \otimes \RT_{\gamma}(H))\circ\Delta_{\alpha\beta,\gamma}\bigr)(v_{\alpha\beta\gamma})
  =\bigl((\Delta_{\alpha,\beta}\otimes H_{\gamma})(\tilde{Q}_{\alpha\beta,\gamma})\bigr)(\tilde{Q}_{\alpha,\beta})_{12\gamma}
  (v_{\alpha}\otimes v_{\beta}\otimes v_{\gamma})$,
 while, by computing the right-hand side,
 $\bigl((\RT_{\alpha}(H)\otimes \Delta_{\beta,\gamma})\circ \Delta_{\alpha,\beta\gamma}\bigr)(v_{\alpha\beta\gamma})
  =\bigl((H_{\alpha}\otimes \Delta_{\beta,\gamma})(\tilde{Q}_{\alpha,\beta\gamma})\bigr)
     (\tilde{Q}_{\beta,\gamma})_{\alpha 23}(v_{\alpha}\otimes v_{\beta}\otimes v_{\gamma})$.
 We only need to show that
 $\bigl((\Delta_{\alpha,\beta}\otimes H_{\gamma})(\tilde{Q}_{\alpha\beta,\gamma})\bigr)(\tilde{Q}_{\alpha,\beta})_{12\gamma}=
  \bigl((H_{\alpha}\otimes \Delta_{\beta,\gamma})(\tilde{Q}_{\alpha,\beta\gamma})\bigr)(\tilde{Q}_{\beta,\gamma})_{\alpha 23}$, 
 or, equivalently,
 \begin{equation}\label{e:tatata}\begin{split}
 \Bigl(\bigl(&\sigma\circ(H_{\beta}\otimes\varphi_{\alpha})\bigr)(R_{\beta,\alpha})\Bigr)_{12\gamma}
 (R_{\alpha,\beta})_{12\gamma}\cdot \\
 & \qquad\cdot\bigl((\Delta_{\alpha,\beta}\otimes H_{\gamma})\circ\sigma\circ
 (H_{\gamma}\otimes\varphi_{\alpha\beta})\bigr)(R_{\gamma,\alpha\beta})
 (\Delta_{\alpha,\beta}\otimes H_{\gamma})(R_{\alpha\beta,\gamma})\\
  & = \Bigl(\bigl(\sigma\circ(H_{\gamma}\otimes\varphi_{\beta})\bigr)(R_{\gamma,\beta})\Bigr)_{\alpha 23}
      (R_{\beta,\gamma})_{\alpha 23}\cdot\\
  & \qquad\cdot\bigl((H_{\alpha}\otimes\Delta_{\beta,\gamma})\circ\sigma\circ
  (H_{\beta\gamma}\otimes\varphi_{\alpha})\bigr)(R_{\beta\gamma,\alpha}) 
  (H_{\alpha}\otimes\Delta_{\beta,\gamma})(R_{\alpha,\beta\gamma})\text{.}
 \end{split}\end{equation}
 Let us set
 $x= \biggl(\Bigl(\bigl(\sigma\circ(\varphi_{\alpha^{-1}}\otimes H_{\alpha})\circ
      \Delta_{\alpha\beta\alpha^{-1},\alpha}\bigr)\otimes H_{\gamma}\Bigr)
      \circ\sigma\circ(H_{\gamma}\otimes\varphi_{\alpha\beta})\biggr)(R_{\gamma,\alpha\beta})$,
 and
 $y= \biggl(\Bigl(H_{\alpha}\otimes\bigl(\sigma\circ(\varphi_{\beta^{-1}}\otimes H_{\beta})\circ
  \Delta_{\beta\gamma\beta^{-1},\beta}\bigr)\Bigr)
  \circ\sigma\circ(H_{\beta\gamma}\otimes\varphi_{\alpha})\biggr)(R_{\beta\gamma,\alpha})$.
 Since, by~\eqref{e:R-a},
 $(R_{\alpha,\beta})_{12\gamma} \bigl((\Delta_{\alpha,\beta}\otimes H_{\gamma})\circ\sigma
           \circ(H_{\gamma}\otimes\varphi_{\alpha\beta})\bigr)(R_{\gamma,\alpha\beta})
           = x(R_{\alpha,\beta})_{12\gamma}$ and
 $(R_{\beta,\gamma})_{\alpha 23}\bigl((H_{\alpha}\otimes\Delta_{\beta,\gamma})
           \circ\sigma\circ(H_{\beta\gamma}\otimes\varphi_{\alpha})\bigr)(R_{\beta\gamma,\alpha})
           = y(R_{\beta,\gamma})_{\alpha 23}$,
 if we substitute these expressions in~\eqref{e:tatata} 
 and we apply axiom~\eqref{e:R-c} to the left-hand side and
 axiom~\eqref{e:R-b} to the right-hand side, we find 
 that~\eqref{e:tatata} can be rewritten as
 \begin{equation*}\begin{split}
  \Bigl(\sigma\circ & (H_{\beta}\otimes\varphi_{\alpha})\bigr)(R_{\beta,\alpha})\Bigr)_{12\gamma}
  x(R_{\alpha,\beta})_{12\gamma}
  \bigl((H_{\alpha}\otimes\varphi_{\beta^{-1}})(R_{\alpha,\beta\alpha\beta^{-1}})\bigr)_{1\beta 3}(R_{\beta,\gamma})_{\alpha 23}\\
  & = \Bigl(\bigl(\sigma\circ(H_{\gamma}\otimes\varphi_{\beta})\bigr)(R_{\gamma,\beta})\Bigr)_{\alpha 23}y(R_{\beta,\gamma})_{\alpha 23}
  (R_{\beta,\gamma})_{\alpha 23}(R_{\alpha,\gamma})_{1\beta 3}(R_{\alpha,\beta})_{12\gamma}\text{.}
 \end{split}\end{equation*}
 Thus, by the Yang-Baxter equation~\eqref{e:YB}, 
 we can rewrite~\eqref{e:tatata} as
 \begin{equation}\label{e:SlavaTebieBojhe}
  \Bigl(\sigma\circ (H_{\beta}\otimes\varphi_{\alpha})\bigr)(R_{\beta,\alpha})\Bigr)_{12\gamma}x=
  \Bigl(\bigl(\sigma\circ(H_{\gamma}\otimes\varphi_{\beta})\bigr)(R_{\gamma,\beta})\Bigr)_{\alpha 23}y\text{.}
 \end{equation}
 
 Given three vector spaces $V_{1}$, $V_{2}$, and $V_{3}$, 
 let us introduce the notation $\sigma_{i,j,k}$
 (with $\{i,j,k\}=\{1,2,3\}\/$) for the permutation
 $V_{1}\otimes V_{2}\otimes V_{3}\to V_{i}\otimes V_{j}\otimes V_{k}$.
  
 If we compute the two factors on the left-hand side 
 in~\eqref{e:SlavaTebieBojhe}, then we have
 \begin{equation*}\begin{split}
  \Bigl(\sigma\circ & (H_{\beta}\otimes\varphi_{\alpha})\bigr)(R_{\beta,\alpha})\Bigr)_{12\gamma}
   = \bigl((\sigma\otimes H_{\gamma})\circ\sigma_{2,3,1}\circ(H_{\gamma}\otimes H_{\beta}\otimes\varphi_{\alpha})\bigr)
      \bigl((R_{\beta,\alpha})_{\gamma 23}\bigr)\\
  & = \bigl((\sigma\otimes H_{\gamma})\circ\sigma_{2,3,1}\circ(H_{\beta}\otimes\varphi_{\beta}\otimes\varphi_{\alpha\beta})\bigr)
      \bigl((R_{\beta,\beta^{-1}\alpha\beta})_{\gamma 23}\bigr)\text{,}
 \end{split}\end{equation*}
 (where in the last passage we used axiom~\eqref{e:R-d}) and
 \begin{equation*}\begin{split}
  x & = \biggl(\Bigl(\bigl(\sigma\circ(\varphi_{\alpha^{-1}}\otimes H_{\alpha})\circ
      \Delta_{\alpha^{-1}\beta\alpha,\alpha}\circ\varphi_{\alpha\beta}\bigr)\otimes H_{\gamma}\Bigr)\circ\sigma\biggr)(R_{\gamma,\alpha\beta})\\
    & = \biggl(\Bigl(\bigl(\sigma\circ(\varphi_{\alpha^{-1}}\otimes H_{\alpha})\circ(\varphi_{\alpha\beta}\otimes\varphi_{\alpha\beta})\bigr)\otimes H_{\gamma}\Bigr)
        \circ\sigma_{2,3,1}\circ(H_{\gamma}\otimes\Delta_{\beta,\beta^{-1}\alpha\beta})\biggr)(R_{\gamma,\alpha\beta})\\
    & = \bigl((\sigma\otimes H_{\gamma})\circ\sigma_{2,3,1}\circ(H_{\beta}\otimes\varphi_{\beta}\otimes\varphi_{\alpha\beta})\bigr)
        \bigl((R_{\gamma,\beta^{-1}\alpha\beta})_{1\beta 3}(R_{\gamma,\beta})_{12(\beta^{-1}\alpha\beta)}\bigr)\text{,}
 \end{split}\end{equation*}
 (where, in the last passage, we used axiom~\eqref{e:R-b}). 
 Similarly, on the right-hand side we have
 \begin{multline*}
  \Bigl(\bigl(\sigma  \circ(H_{\gamma}\otimes\varphi_{\beta})\bigr)(R_{\gamma,\beta})\Bigr)_{\alpha 23}
   = \bigl((H_{\alpha}\otimes\sigma)\circ\sigma_{3,1,2}\circ(H_{\gamma}\otimes\varphi_{\beta}\otimes\varphi_{\alpha\beta})\bigr)
    \bigl((R_{\gamma,\beta})_{12(\beta^{-1}\alpha\beta)}\bigr)
 \end{multline*}
 and
 \begin{equation*}\begin{split}
   y & = \biggl((H_{\alpha}\otimes\sigma)\circ\Bigl(\varphi_{\alpha}\otimes\bigl((\varphi_{\beta^{-1}}\otimes H_{\beta})\circ
        \Delta_{\beta\gamma\beta^{-1},\beta}\bigr)\Bigr)\circ\sigma\biggr)(R_{\beta\gamma,\alpha})\\
    & = \Bigl((H_{\alpha}\otimes\sigma)\circ(\varphi_{\alpha\beta}\otimes H_{\gamma}\otimes\varphi_{\beta})
        \circ\bigl(\varphi_{\beta^{-1}}\otimes(\Delta_{\gamma,\beta}\circ\varphi_{\beta^{-1}})\bigr)\circ\sigma\Bigr)(R_{\beta\gamma,\alpha})\\
 \intertext{(by axiom~\eqref{e:R-d})}
    & = \bigl((H_{\alpha}\otimes\sigma)\circ(\varphi_{\alpha\beta}\otimes H_{\gamma}\otimes\varphi_{\beta})
        \circ(H_{\beta^{-1}\alpha\beta}\otimes\Delta_{\gamma,\beta})\circ\sigma\bigr)(R_{\gamma\beta,\beta\alpha\beta^{-1}})\\
    & = \bigl((H_{\alpha}\otimes\sigma)\circ(\varphi_{\alpha\beta}\otimes H_{\gamma}\otimes\varphi_{\beta})
        \circ\sigma_{3,1,2}\circ(\Delta_{\gamma,\beta}\otimes H_{\beta^{-1}\alpha\beta})
        \bigr)(R_{\gamma\beta,\beta\alpha\beta^{-1}})\\
    & = \bigl((H_{\alpha}\otimes\sigma)\circ\sigma_{3,1,2}\circ(H_{\gamma}\otimes\varphi_{\beta}\otimes\varphi_{\alpha\beta})\bigr)
        \Bigl(\bigl((H_{\gamma}\otimes\varphi_{\beta^{-1}})(R_{\gamma,\alpha})\bigr)_{1\beta 3}
        (R_{\beta,\beta^{-1}\alpha\beta})_{\gamma 23}\Bigr)
 \end{split}\end{equation*}
 (where in the last passage we used axiom~\eqref{e:R-c}).
 
 We observe that the application
  $(H_{\alpha}\otimes\sigma)\circ\sigma_{3,1,2}\circ(H_{\alpha}\otimes\varphi_{\alpha}\otimes\varphi_{\alpha\beta})$
 is bijective. 
 So, we can rewrite~\eqref{e:SlavaTebieBojhe} in the form
 \begin{multline*}
  (R_{\beta,\beta^{-1}\alpha\beta})_{\gamma 23}(R_{\gamma,\beta^{-1}\alpha\beta})_{1\beta 3}(R_{\gamma,\beta})_{12(\beta^{-1}\alpha\beta)}\\
   =(R_{\gamma,\beta})_{12(\beta^{-1}\alpha\beta)}\bigl((H_{\gamma}\otimes\varphi_{\beta^{-1}})(R_{\gamma,\alpha})\bigr)_{1\beta 3}(R_{\beta,\beta^{-1}\alpha\beta})_{\gamma 23}
 \end{multline*}
 and this last formula is true by the Yang-Baxter equation~\eqref{e:YB}.
 \end{sentence}
\end{proof}

\begin{cor}
 Let $H$ be a finite-type \Tcoalg. We obtain a ribbon \Tcoalg 
 $\RT\bigl(D(H)\bigr)$.
\end{cor}

\bibliography{AlgebraicDouble}

\end{document}